\documentclass[11pt,a4paper]{amsart}
\usepackage[english]{babel}
\usepackage{graphicx}
\usepackage{framed}
\usepackage[normalem]{ulem}
\usepackage{amsmath}
\usepackage{amsthm}
\usepackage{amssymb}
\usepackage{amsfonts}
\usepackage[foot]{amsaddr}
\usepackage{mathrsfs,amsmath}
\usepackage{enumerate}
\usepackage[utf8]{inputenc}
\usepackage[top=1 in,bottom=1in, left=1 in, right=1 in]{geometry}
\usepackage{mathtools}
\usepackage{tikz-cd}
\usepackage[font=small,labelfont=bf]{caption} 
\usepackage[subrefformat=parens]{subcaption}
\usepackage{cite}
\usepackage{abstract}

\usetikzlibrary{decorations.markings, arrows.meta}

\tikzset{
  marrow/.style={decoration={markings,mark=at position 0.5 with {\arrow{#1}}}, postaction=decorate}
}

\usepackage{enumitem}
\setlist{  
  listparindent=\parindent,
  parsep=0pt,
}

\usepackage{bbm}

\usepackage{stackrel}

\newcommand{\supp}{{\rm supp}}

\newcommand{\R}{\mathbb{R}}
\newcommand{\N}{\mathbb{N}}

\newcommand{\Z}{\mathbb{Z}}

\newcommand{\eps}{\varepsilon}

\newtheorem{thm}{Theorem}[section]
\newtheorem{cor}[thm]{Corollary}
\newtheorem{lem}[thm]{Lemma}
\newtheorem{que}[thm]{Question}
\newtheorem{prop}[thm]{Proposition}

\theoremstyle{definition}

\theoremstyle{remark}
\newtheorem{ex}[thm]{Example}

\newtheorem{rmk}[thm]{Remark}

\newcommand{\bproof}{\noindent{\textit{Proof. }}}
\newcommand{\eproof}{\hfill\qed}

\newcommand{\Cont}{{\rm Cont_0}}
\newcommand{\tildecont}{{\rm \widetilde{Cont}_0}}

\DeclareFontFamily{U}{mathb}{\hyphenchar\font45}
\DeclareFontShape{U}{mathb}{m}{n}{
	<-6> mathb5 <6-7> mathb6 <7-8> mathb7
	<8-9> mathb8 <9-10> mathb9
	<10-12> mathb10 <12-> mathb12
}{}
\DeclareSymbolFont{mathb}{U}{mathb}{m}{n}
\DeclareMathSymbol{\lprec}{\mathrel}{mathb}{"CE}
\DeclareMathSymbol{\gprec}{\mathrel}{mathb}{"CF}

\begin{document}

\title{A new metric on the contactomorphism group of orderable contact manifolds}
\author{Lukas Nakamura}
\maketitle

\begin{abstract}
We introduce a pseudo-metric on the contactomorphism group of any contact manifold $(M,\xi)$ with a cooriented contact structure $\xi$. It is the contact analogue of a corresponding semi-norm in Hofer's geometry, and on certain classes of contact manifolds, its lift to the universal cover can be viewed as a continuous version of the integer valued bi-invariant metric introduced by Fraser, Polterovich, and Rosen. We show that it is non-degenerate if and only if $(M,\xi)$ is strongly orderable and that its metric topology agrees with the interval topology introduced by Chernov and Nemirovski. In particular, the interval topology is Hausdorff whenever it is non-trivial, which answers a question of Chernov and Nemirovski. We discuss analogous results for isotopy classes of Legendrians and their universal covers.
\end{abstract}

\section{Introduction}

Let $(M, \xi)$ be a connected, cooriented contact manifold of dimension $2n+1$, and let $\Cont(M,\xi)$ denote identity component of the space of compactly supported contactomorphisms on $M$. Let $\alpha$ be a contact form for $\xi$ so that the induced coorientation agrees with the given one. The contact form $\alpha$ induces a one-to-one correspondence between compactly supported contact vector fields on $(M,\xi)$ and functions $H:M \to \R$ by associating to a contact vector field $X$ the function $H = \alpha(X)$, called the \emph{Hamiltonian}. Conversely, to a (time-dependent) Hamiltonian $H_t$ we can associate a (time-dependent) contact vector field $X^H_t$, and by integrating $X^H_t$ we obtain a contact isotopy $\phi^H_t$ starting at the identity.\\

With the goal of studying the geometry of $\Cont(M,\xi)$, Eliashberg and Polterovich \cite{ep00} introduced the binary relation $\leq$ on the universal cover $\tildecont(M,\xi)$ of the contactomorphism group $\Cont(M,\xi)$, where $f \leq g$ if and only if there exists a \emph{non-negative} contact isotopy from $f$ to $g$ in the sense that the associated Hamiltonian is non-negative. It is clear that $\leq$ is transitive and reflexive, and that it is antisymmetric if and only if there does not exists any contractible non-trivial non-negative loop of contactomorphisms on $(M,\xi)$. In the case that any contractible non-negative loop of contactomorphisms is the constant loop, $(M,\xi)$ is called \emph{orderable}. If $M$ is compact, one similarly defines the quasi-order\footnote{Note that $\lprec$ is not reflexive unless there exists a positive loop of contactomorphisms. We will still refer to $\lprec$ as a quasi-order.} $\lprec$ on $\tildecont(M,\xi)$ via \emph{positive} contact isotopies generated by a \emph{positive} Hamiltonian. Eliashberg and Polterovich showed that $(M,\xi)$ is orderable if and only if there does not exist a contractible positive loop of contactomorphisms in which case also $\lprec$ defines a partial order on $\tildecont(M,\xi)$.

Using this partial order, they defined \emph{the relative growth} of two elements $f,g \in \tildecont(M,\xi)$ with $f \gprec Id_M$  as 
\begin{equation}
\gamma(f,g) \coloneqq \lim\limits_{k \to \infty} \frac{\inf\{p \in \Z|\, f^p \geq g^k\}}{k}
\end{equation} 
and showed that the limit exists in $\R$ whenever $(M,\xi)$ is orderable. \\

These constructions extend naturally to $\Cont(M,\xi)$ to give quasi-orders $\leq$ and $\lprec$ and the relative growth $\gamma(\phi,\psi)$ of two contactomorphisms $\phi,\psi \in \Cont(M,\xi)$ with $\phi \gprec Id_M$. On $\Cont(M,\xi)$, both $\leq$ and $\lprec$ are anti-symmetric if and only if there does not exists \emph{any} positive loop of contactomorphisms on $(M,\xi)$ in which case $(M,\xi)$ is called \emph{strongly orderable}, and in this case $\gamma(\phi,\psi) \in \R$ always exists. Following \cite{ep00}, (strong) orderability or non-orderability was shown for many contact manifolds (see for example \cite{af12}, \cite{afm15}, \cite{am18}, \cite{bhu01}, \cite{bz15}, \cite{cfp17}, \cite{cn10b}, \cite{cn16}, \cite{ekp06}, \cite{wei13}).\\

Now fix a closed Legendrian $L \subseteq M$, and recall that any isotopy of $L$ through Legendrians is induced by an ambient contact isotopy. As was already noted by Eliashberg and Polterovich, this implies that $\leq$ also induces a quasi-order on the Legendrian isotopy class $\mathcal{L}(L)$ of $L$ and its universal cover $\widetilde{\mathcal{L}}(L)$, and if $L$ is closed, also $\lprec$ induces a quasi-order. We will denote the induced quasi-orders by the by the same symbols $\leq$ and $\lprec$. $L$ is called (strongly) orderable if $\leq$ is a partial order on $\widetilde{\mathcal{L}}(L)$ (resp. $\mathcal{L}(L)$). As before, in the compact case $L$ is (strongly) orderable if and only if there does not exist a positive loop of Legendrians in $\widetilde{\mathcal{L}}(L)$ (resp. $\mathcal{L}(L)$). These quasi-orders were first studied in detail in \cite{cfp17} where among other results it was shown that the fibre of the unit cotangent bundle $ST^*N$ of any manifold $N$ that is smoothly covered by $\R^n$ is strongly orderable. In \cite{cn10b} this result was extended to the case that the universal cover of $N$ is an open manifold, and in \cite{cn16} it was shown for arbitrary $N$ that the fibre in $ST^*N$ is orderable.\\

In analogy to the Hofer metric on the group of Hamiltonian diffeomorphisms of a symplectic manifold, one can associate to a contact isotopy an energy which is given by the $L^{1,\infty}$-norm of the corresponding Hamiltonian:
\begin{equation}
\Vert \phi^H_t \Vert_\alpha := \Vert H_t \Vert \coloneqq \int\limits_0^1 \max\limits_M |H_t| dt.
\end{equation}

It follows from the behaviour of the Hamiltonians under left-/right-composition and inversion of contact isotopies that the induced \emph{Shelukhin-Hofer metric} 
\begin{equation}
d^\alpha_{SH}(\phi_0,\phi_1) \coloneqq \inf \{\Vert H_t \Vert |\, H_t:M \to \R,\, \phi^H_1 \phi_0 = \phi_1 \}
\end{equation}
defines a right-invariant pseudo-metric on $\Cont(M,\xi)$ which behaves naturally under left-composition:
\begin{equation}
d^\alpha_{SH}(\psi \phi_0,\psi \phi_1) = d^{\psi^*\alpha}_{SH}(\phi_0,\phi_1).
\end{equation}

Furthermore, it follows easily that $d^\alpha_{SH}$ and $d^{\alpha'}_{SH}$ are equivalent for any two choices of contact forms $\alpha$ and $\alpha'$ for $\xi$ which agree outside of a compact subset.

In \cite{she17}, Shelukhin showed that $d^\alpha_{SH}$ is, in fact, non-degenerate on any contact manifold. Consequently, it defines a genuine metric on $\Cont(M,\xi)$.\\

The energy of a contact isotopy also naturally induces a right-invariant and left-natural pseudo-metric $\widetilde{d}^\alpha_{SH}$ on the universal cover $\tildecont(M,\xi)$ of $\Cont(M,\xi)$, but as was pointed out by Shelukhin \cite{she17}, it is only known that $\widetilde{d}^\alpha_{SH}(f,g) > 0$ for $f,g \in \tildecont(M,\xi)$ if $\pi(f) \neq \pi(g)$ where $\pi: \tildecont(M,\xi) \to \Cont(M,\xi)$ denotes the universal covering map.\\

$d^\alpha_{SH}$ (and $\widetilde{d}^\alpha_{SH}$) induce pseudo-metrics $d^\alpha_{SCH}$ (and $\widetilde{d}^\alpha_{SCH}$) on orbits of parametrized or unparametrized subsets $A$ of $M$ (and their universal cover) under the action of $\Cont(M,\xi)$ (and $\tildecont(M,\xi)$), which we will refer to as the \emph{Shelukhin-Chekanov-Hofer metric}. $d^\alpha_{SCH}$ was studied for closed unparametrized submanifolds $L$ by Rosen and Zhang in \cite{rz18} where they observed the following. The pseudo-metric $d^\alpha_{SCH}$ can only be non-degenerate if $L$ is coisotropic, and if $\dim(L) < n$, where $\dim(M) = 2n+1$, then $d^\alpha_{SCH}$ vanishes identically. In the case that $\dim(L) = n$, $d^\alpha_{SCH}$ vanishes identically unless $L$ has a Legendrian component, and if $L$ is a connected Legendrian, then $d^\alpha_{SCH}$ either vanishes identically or is non-degenerate. For disjoint unions of connected, parametrized non-Legendrian submanifolds of dimension $n$, Dimitroglou Rizell and Sullivan \cite{drs24b} proved that $d^\alpha_{SCH}$ vanishes identically as well.\\

Motivated by the Alexandrov topology in Lorentzian geometry, Chernov and Nemirovski \cite{cn20} used the quasi-order $\lprec$ to define the \emph{interval topologies} $\widetilde{\mathcal{O}}_\lprec$ and $\mathcal{O}_\lprec$ on $\tildecont(M,\xi)$ and $\Cont(M,\xi)$ respectively, which are generated by intervals
\begin{equation}\label{eq:defintion of intervals}
I_{f_0, f_1}^\lprec = \{f \in \tildecont(M,\xi)|\, f_0 \lprec f \lprec f_1\} 
\end{equation}
in $\tildecont(M,\xi)$ with respect to $\lprec$, and similarly for $\mathcal{O}_\lprec$. In fact, it is straightforward to see that the intervals $I^\lprec_{f^\alpha_{-\eps} f,f^\alpha_{\eps} f}$ and $I^\lprec_{\phi^\alpha_{-\eps} \phi,\phi^\alpha_{\eps}\phi}$, $\eps > 0, f \in \tildecont(M,\xi), \phi \in \Cont(M,\xi)$, form a basis of $\widetilde{\mathcal{O}}_\lprec$ and $\mathcal{O}_\lprec$, respectively, where $\phi^\alpha_t$ denotes the Reeb flow and $f^\alpha_t \coloneqq [s \mapsto \phi^\alpha_{st}]$ denotes its lift to $\tildecont(M,\xi)$.

It follows easily that the interval topologies are trivial if there exists a (contractible) positive loop of contactomorphisms. In the absence of positive loops, Chernov and Nemirovski asked whether the interval topologies are necessarily Hausdorff.

\begin{que}\label{que:orderable=hausdorff}\hspace{-1mm}\textnormal{\cite{cn20}}\,
	 If $(M,\xi)$ is (strongly) orderable, is $\widetilde{\mathcal{O}}_\lprec$ ($\mathcal{O}_\lprec$) Hausdorff?
\end{que}

They proved that this is indeed the case on spherical cotangent bundles of closed manifolds which can be smoothly covered by an open subset of $\R^n$.\\

Hedicke \cite{hed21} showed that the interval topology  is coarser than the topology $\mathcal{O}_{d^\alpha_{SH}}$ induced by $d_{SH}^\alpha$, i.e. $\mathcal{O}_\lprec \subseteq \mathcal{O}_{d^\alpha_{SH}}$, which leads to the following natural question.

\begin{que}\hspace{-1mm}\textnormal{(Question 2 in \textnormal{\cite{hed21}})}\,
	Are $\mathcal{O}_\lprec$ and $\mathcal{O}_{d^\alpha_{SH}}$ equal?\\
\end{que}

Using the quasi-order $\lprec$ on $\mathcal{L}(L)$ and $\widetilde{\mathcal{L}}(L)$, Chernov and Nemirovski \cite{cn20} again defined intervals $I^\lprec_{L_0,L_1}$ and interval topologies $\mathcal{O}_\lprec^L$ and $\widetilde{\mathcal{O}}_\lprec^L$ for a closed Legendrian $L \subseteq M$. As above, these can only be non-trivial if $L$ is (strongly) orderable, and intervals whose endpoints are related via the Reeb flow form a basis. Hedicke \cite{hed21} proved also in this case that $\mathcal{O}^L_\lprec \subseteq \mathcal{O}_{d^\alpha_{CSH}}$.\\

Further inspired by Lorentzian geometry, Hedicke \cite{hed21} recently introduced the Lorentzian distance function
\begin{equation}
\begin{split}
\tau_\alpha: \Cont(M,\xi) \times \Cont(M,\xi) &\to [0, \infty]\\
(\phi_0,\phi_1) &\mapsto \begin{cases} \sup \left\{ \int\limits_0^1 \min\limits_M H_t dt|\,  \phi^H_1 \phi_0 = \phi_1, H \geq 0 \right\} \quad &\text{ if } \phi_1 \geq \phi_0, \\ 0 \quad &\text{ otherwise.} \end{cases}
\end{split}
\end{equation}
and showed that if $(M,\xi)$ is strongly orderable, then $\tau_\alpha < \infty$, and that $\tau_\alpha(\phi_0,\phi_1) > 0$ if and only if $\phi_0 \lprec \phi_1$. Furthermore, he proved that $\tau_\alpha$ is continuous with respect to $\mathcal{O}_\lprec$ and thus also with respect to $d_{SH}^\alpha$. \\

As explained in \cite{hed21}, $\tau_\alpha$ can also be defined on $\mathcal{L}(L)$ for any closed Legendrian $L$, and it satisfies similar properties as on $\Cont(M,\xi)$.\\

The definition of $\tau_\alpha$ extends naturally also to $\tildecont(M,\xi)$ and $\widetilde{\mathcal{L}}(L)$, and it is easy to check Hedicke's results mentioned above still apply with the same proofs.\\

The goal of this article is to further analyse the above questions. Similarly to the semi-norms in symplectic geometry which measure only the positive or negative part of a normalized Hamiltonian, we introduce below ``distances" $\delta^\pm_\alpha$ on $\Cont(M,\xi)$ which are bounded by the Shelukhin-Hofer metric and can be used to express the Lorentzian distance function. It turns out that their maximum $d_\alpha \coloneqq \max\{\delta^\pm_\alpha,0\}$ defines a pseudo-metric on $\Cont(M,\xi)$. In fact, we prove the following, which answers Question \ref{que:orderable=hausdorff} above (see Section \ref{sec:cont(M,xi)} below).

\begin{thm}\label{thm:main result}
	$d_\alpha$ either vanishes identically or is  non-degenerate, and it satisfies $d_\alpha \leq d^\alpha_{SH}$. If $M$ is compact, then the topology on $\Cont(M,\xi)$ induced by $d_\alpha$ agrees with the interval topology $\mathcal{O}_\lprec$, and $d_\alpha$ vanishes identically if and only if there exists a positive loop of contactomorphisms on $M$. \eproof\\
\end{thm}

We also introduce versions of $\delta^\pm_\alpha$ and $d_\alpha$ on $\tildecont(M,\xi)$ and spaces of parametrized or unparametrized subspaces and their universal covers and discuss similar results for these spaces, but similarly to the Shelukhin-Hofer metric, our methods do not suffice to prove genuine non-degeneracy on the universal covers.\\

\begin{rmk}
	Allais and Arlove \cite{aa23} independently also studied the pseudo-metrics considered here and obtained proofs of our main results. After prepublication of the first version of this article, it came to our attention that the definition of the pseudo-norms associated to $\widetilde{\delta}^\pm_\alpha$ and $\widetilde{d}_\alpha$ on $\tildecont(M,\xi)$ together with Example \ref{ex:deltapm on ST*Tn} below were already introduced in Section 5 of Arlove's PhD thesis \cite{arl21}.\\
\end{rmk}

\textbf{Acknowledgements:} I am grateful to Georgios Dimitroglou Rizell for introducing me to Question \ref{que:orderable=hausdorff} which was the starting point for this project, and to Stefan Nemirovski for explaining to me part of his work on orderability. I would also like to thank my advisor Tobias Ekholm for valuable discussions and for careful reading of a preliminary version of this article. I want to thank Simon Allais and Pierre-Alexandre Arlove for directing my attention to Arlove's PhD thesis. I am also grateful to the anonymous referee for helpful comments on an earlier version of this article.

This work is partially supported by the Knut and Alice Wallenberg Foundation, grant KAW2020.0307.
\section{Main results}

\subsection{Distances on the universal cover $\tildecont(M,\xi)$ of $\Cont(M,\xi)$}\hfill 

For any two $f_0, f_1 \in \tildecont(M,\xi)$, we define 
\begin{equation}
\begin{split}
\widetilde{\delta}^-_\alpha(f_0,f_1) &\coloneqq \inf \{\eps \in \R|\, \exists H_t:M \to [-\eps,\infty): [\phi^{H_t}_t] \cdot f_0 = f_1\},\\
\widetilde{\delta}^+_\alpha(f_0,f_1) &\coloneqq \inf \{\eps \in \R|\, \exists H_t:M \to (\infty, \eps]: [\phi^{H_t}_t] \cdot f_0 = f_1\}.
\end{split}
\end{equation}

Recall that we write $f \lprec g$ for $f,g \in \tildecont(M,\xi)$ if there exists a positive contact isotopy from $f$ to $g$. Using the class $f^\alpha_t = [s \mapsto \phi^\alpha_{st}, s \in [0,1]] \in \tildecont(M,\xi)$ of the time-$t$ Reeb flow $\phi^\alpha_t$ and the cones 
\begin{equation}
I^\pm_f \coloneqq \left\{g \in \tildecont(M,\xi)|\, f^{\pm1} \lprec g^{\pm1} \right\}
\end{equation} in the case that $M$ is compact, $\widetilde{\delta}^\pm_\alpha$ can also be expressed as 
\begin{equation}\label{eq:deltapm via cones}
\widetilde{\delta}^\pm_\alpha(f_0,f_1) = \inf \left\{\eps \in \R|\, f_1 \in I^\mp_{f^\alpha_{\pm \eps} f_0} \right\},
\end{equation}
since any contact isotopy $\phi_t$ can be written as $\phi^\alpha_{-\eps t} \circ (\phi^\alpha_{\eps t} \circ \phi_t) = \phi^\alpha_{\eps t} \circ (\phi^\alpha_{-\eps t} \circ \phi_t)$ and the Reeb flow is a strict contact isotopy. Then the Hamiltonian associated to $\phi_t$ is $\geq -\eps$ if and only if the Hamiltonian associated to $\phi^\alpha_{\eps t} \circ \phi_t$ is $\geq 0$, and it is $\leq \eps$ if and only if the Hamiltonian associated to $\phi^\alpha_{-\eps t} \circ \phi_t$ is $\leq 0$. 

This together with the observation that for any $\eta > 0$ and $f_1 \in I^\mp_{f^\alpha_{\pm \eps} f_0}$, we have $I^\lprec_{f^\alpha_{-\eta} f_1,f^\alpha_\eta f_1} = I^-_{f^\alpha_\eta f_1} \cap I^+_{f^\alpha_{-\eta} f_1} \subseteq I^\mp_{f^\alpha_{\pm (\eps + \eta)}f_0}$ (see (\ref{eq:defintion of intervals}) above for the definition of $I^\lprec$) shows that $\widetilde{\delta}^\pm_\alpha$ is continuous with respect to the interval topology. Note also that $\widetilde{\delta}^-_\alpha(f_0,f_1) < 0$ if and only if $f_1 \gprec f_0$.\\

Let $\pi:\tildecont(M,\xi) \to \Cont(M,\xi)$ denote the natural projection. $\widetilde{\delta}^+_\alpha$ and $\widetilde{\delta}^-_\alpha$ satisfy the following properties.

\begin{prop}\label{prop:properties of deltapm} Let $f_0, f_1, f_2 \in \tildecont(M,\xi)$. Then
	
	(i) (triangle inequality) $\widetilde{\delta}^\pm_\alpha(f_0, f_2) \leq \widetilde{\delta}^\pm_\alpha(f_0,f_1) + \widetilde{\delta}^\pm_\alpha(f_1,f_2)$.
	
	(ii) (symmetry) $\widetilde{\delta}^+_\alpha(f_0, f_1) = \widetilde{\delta}^-_\alpha(f_1, f_0)$.
	
	(iii) (left naturality) $\widetilde{\delta}^\pm_\alpha(f_2 f_0,f_2 f_1) = \widetilde{\delta}^\pm_{(\pi(f_2))^* \alpha}(f_0,f_1).$
	
	(iv) (inverse) $\widetilde{\delta}^+_\alpha(id_M, f_0^{-1}) = \widetilde{\delta}^-_\alpha(id_M, f_0,)$.
	
	(v) (right invariance) $\widetilde{\delta}^\pm_\alpha(f_0 \circ f_2,f_1 \circ f_2) = \widetilde{\delta}^\pm_\alpha(f_0, f_1)$.
	
	(vi) (orderability) Assume that $M$ is compact. If $(M,\xi)$ is not orderable, then $\widetilde{\delta}^\pm_\alpha \equiv - \infty$. If $(M,\xi)$ is orderable, then $\widetilde{\delta}^\pm_\alpha > - \infty$, and at least one of $\widetilde{\delta}^+_\alpha(f_0,f_1)$ and $\widetilde{\delta}^-_\alpha(f_0,f_1)$ is non-negative, and if one is negative, then the other one is positive.\\
\end{prop}

\bproof 
The properties $(ii) - (v)$ follow similarly to the corresponding properties of the Shelukhin-Hofer metric $d^\alpha_{SH}$ from the behavior of the Hamiltonians (see the proof of Theorem A in \cite{she17}). $(i)$ is a consequence of Proposition \ref{prop:delta as integral} below. To see $(vi)$, first note that if there exists a contractible positive loop of contactomorphisms, then for any $f \in \tildecont(M,\xi)$, $I^\pm_f = \tildecont(M,\xi)$, so $\widetilde{\delta}^\pm_\alpha = - \infty$. Conversely, assume that $\widetilde{\delta}^\pm_\alpha(f_0,f_1) = -\infty$. Then by the triangle inequality, also $\widetilde{\delta}^\pm_\alpha(f_1,f_1) = -\infty$. If $\widetilde{\delta}^+_\alpha(f_0,f_1) < 0$ and $\widetilde{\delta}^-_\alpha(f_0,f_1) = \widetilde{\delta}^+_\alpha(f_1,f_0) \leq 0$ or $\widetilde{\delta}^+_\alpha(f_0,f_1) \leq 0$ and $\widetilde{\delta}^-_\alpha(f_0,f_1) = \widetilde{\delta}^+_\alpha(f_1,f_0) < 0$, then again by the triangle inequality, $\widetilde{\delta}^\pm_\alpha(f_1,f_1) < 0$. In either case $f_1 \gprec f_1$, so $(M,\xi)$ is not orderable.
\eproof\\

Now we define
\begin{equation}
\widetilde{d}_\alpha(f_0,f_1) \coloneqq \max \{\widetilde{\delta}^-_\alpha(f_0,f_1),\widetilde{\delta}^+_\alpha(f_0,f_1),0\}.
\end{equation}\\

\begin{rmk}
In the case that the Reeb flow of $(M,\alpha)$ is $1$-periodic, this definition was previously considered in the form of the metric balls of $\widetilde{d}_\alpha$ in \cite{fpr18} where it was shown that the boundary of the balls 
\begin{equation}
\mathcal{B}(k) \coloneqq \{f \in \tildecont(M,\xi)|\, f^\alpha_{-k} \leq f \leq f^{\alpha}_k\}, \quad k \in \N,
\end{equation} with respect to the interval topology $\widetilde{\mathcal{O}}_\lprec$ only contains lifts of contactomorphisms with $k$- or $-k$-translated points.\\
\end{rmk}

In the case that $M$ is compact, it follows from \eqref{eq:deltapm via cones} that $\widetilde{d}_\alpha(f_0,f_1)$ is also equal to 
\begin{equation}
\widetilde{d}_\alpha(f_0,f_1) = \inf \{\eps > 0|\, f_1 \in I^\lprec_{f^\alpha_{-\eps}f_0, f^\alpha_\eps f_0} \},
\end{equation} 
or equivalently,
\begin{equation}\label{eq:metric ball is interval}
\{f_1 \in \tildecont(M,\xi)|\, \widetilde{d}_\alpha(f,f_1) < \eps \} = I^\lprec_{f^\alpha_{-\eps}f, f^\alpha_\eps f}.
\end{equation}

With Proposition \ref{prop:properties of deltapm}, it is immediate that $\widetilde{d}_\alpha$ is a pseudo-metric, and by (\ref{eq:metric ball is interval}), its metric topology agrees with the interval topology $\widetilde{\mathcal{O}}_\lprec$ since the intervals $I^\lprec_{f^\alpha_{-\eps}f, f^\alpha_\eps f}$ form a basis. 

In fact, the following holds.

\begin{prop}\label{prop:properties of d_alpha} Let $f_0, f_1, f_2 \in \tildecont(M,\xi)$. Then
	
	(i) (triangle inequality) $\widetilde{d}_\alpha(f_0, f_2) \leq \widetilde{d}_\alpha(f_0,f_1) + \widetilde{d}_\alpha(f_1,f_2)$.
	
	(ii) (symmetry) $\widetilde{d}_\alpha(f_0, f_1) = \widetilde{d}_\alpha(f_1, f_0)$.
	
	(iii) (left naturality) $\widetilde{d}_\alpha(f_2 f_0,f_2 f_1) = d_{(\pi(f_2))^* \alpha}(f_0,f_1).$
	
	(iv) (inverse) $\widetilde{d}_\alpha(f_0^{-1}, id_M) = \widetilde{d}_\alpha(f_0,id_M)$.
	
	(v) (right invariance) $\widetilde{d}_\alpha(f_0 \circ f_2,f_1 \circ f_2) = \widetilde{d}_\alpha(f_0, f_1)$.
	
	(vi) (topology) If $M$ is compact, then the metric topology of $\widetilde{d}_\alpha$ agrees with the interval topology $\widetilde{\mathcal{O}}_\lprec$.
\end{prop}

\bproof 
Immediate from Proposition \ref{prop:properties of deltapm} and the above discussion.
\eproof\\

The properties in Proposition \ref{prop:properties of d_alpha} together with the following straightforward lemma imply that $(\tildecont(M,\xi),\widetilde{d}_\alpha)$ is a topological group.

\begin{lem}\label{lem:differrent alpha}
	Let $\alpha' = \lambda \alpha$ be another contact form for $\xi$ where $\lambda:M \to \R_{> 0}$ is a strictly positive function that is equal to the constant function $1$ outside of a compact set. Then
	\begin{equation}
	\underset{M}{\min}(\lambda)\, \widetilde{d}_\alpha \leq \widetilde{d}_{\alpha'} \leq \underset{M}{\max}(\lambda)\, \widetilde{d}_\alpha.
	\end{equation}
	In particular, $\widetilde{d}_\alpha$ and $\widetilde{d}_{\alpha'}$ induce the same topology on $\tildecont(M)$.
\end{lem}

\bproof
Clear.
\eproof\\


It is also easy to see that $\widetilde{d}_\alpha$ is always bounded by the Shelukhin-Hofer metric $\widetilde{d}^\alpha_{SH}$. In fact, if we denote by $\widetilde{\delta}_{SH}^+(f_0,f_1)$, resp. $\widetilde{\delta}_{SH}^-(f_0,f_1)$, the infimum of $\Vert H_t \Vert^+$, resp. $\Vert H_t \Vert^-$, over all Hamiltonians $H_t$ whose associated contact isotopy $\phi^H_t$ satisfies $[\phi^H_t] \cdot f_0 = f_1$, where $\Vert H_t \Vert^\pm \coloneqq \int\limits_0^1 \underset{M}{\max}\{\pm H_t\} dt$, we have the following.

\begin{prop}\label{prop:delta as integral}
	For any $f_0, f_1 \in \tildecont(M,\alpha)$, 
	\begin{equation}\label{eq:delta-sh=delta-alpha}
	\widetilde{\delta}^{\pm}_{SH} = \widetilde{\delta}^{\pm}_\alpha.
	\end{equation}
	In particular,
	\begin{equation}
	\widetilde{d}_\alpha(f_0,f_1) = \max \left\{\widetilde{\delta}_{SH}^-(f_0,f_1),\widetilde{\delta}_{SH}^+(f_0,f_1), 0 \right\} \leq d_{SH}(f_0,f_1).
	\end{equation}
\end{prop} 

\bproof
Note that if $M$ is non-compact, $\underset{M}{\max}\{\pm H_t\} \geq 0$ for any compactly supported Hamiltonian $H_t$, in which case Proposition \ref{prop:delta as integral} follows directly from Lemma \ref{lem:reparam of hamiltonian to make max const} below applied to $\underset{M}{\max}\{\pm H_t\}$. In the compact case, the \eqref{eq:delta-sh=delta-alpha} is just  (\ref{eq:tau,delta in terms of intergrals}) in Section \ref{sec:hedickes lorentzian distance function} below.
\eproof\\

Next we discuss the non-degeneracy of $\widetilde{d}_\alpha$. Part $(vi)$ of Proposition \ref{prop:properties of deltapm} shows that $\widetilde{d}_\alpha$ vanishes identically if $M$ is compact and $(M,\xi)$ is not orderable. It turns out that the converse is true as well, and that $\widetilde{d}_\alpha$ is always non-degenerate in a weaker sense unless it vanishes identically. We say that a pseudo-metric $d$ on $\tildecont(M,\xi)$ is $\pi$-non-degenerate if for any $f_0,f_1 \in \tildecont(M,\xi)$ with $\pi(f_0) \neq \pi(f_1)$, we have that $d(f_0,f_1) > 0$. Similarly, we say that a topology $\mathcal{O}$ on $\tildecont(M,\xi)$ is $\pi$-Hausdorff if any $f_0,f_1 \in \tildecont(M,\xi)$ with $\pi(f_0) \neq \pi(f_1)$ have disjoint open neighbourhoods, where $\pi: \tildecont(M,\xi) \to \Cont(M,\xi)$ is the covering map.  Clearly, these notions are weaker than Hausdorffness and, respectively, non-degeneracy in general. 

\begin{thm}\label{thm:non-deg of d_alpha} 
	$\widetilde{d}_\alpha$ either vanishes identically or is  $\pi$-non-degenerate. If $M$ is compact, then $\widetilde{d}_\alpha$ vanishes identically if and only if there exists a contractible positive loop of contactomorphisms on $M$.\\
\end{thm}

The proof of Theorem \ref{thm:non-deg of d_alpha} is given in Section \ref{sec:proofs of main results}.\\


Since $\widetilde{\mathcal{O}}_\lprec$ is the metric topology of $\widetilde{d}_\alpha$, it is clear that $\widetilde{\mathcal{O}}_\lprec$ is ($\pi$-)Hausdorff if $\widetilde{d}_\alpha$ is ($\pi$-)non-degenerate. An immediate consequence is the following.

\begin{cor}
	Let $(M,\xi)$ be a closed contact manifold. Then the interval topology $\widetilde{\mathcal{O}}_\lprec$ on $\tildecont(M,\xi)$ is $\pi$-Hausdorff if and only if there does not exist a contractible positive loop of contactomorphisms on $(M,\xi)$.\\
\end{cor}

In the case that $\widetilde{d}_\alpha$ is $\pi$-non-degenerate, we now have the two natural metrics $\widetilde{d}_\alpha \leq \widetilde{d}^\alpha_{SH}$ on $\tildecont(M,\xi)$. 

\begin{que}
	Is $\widetilde{d}_\alpha$ equivalent to $\widetilde{d}^\alpha_{SH}$ if it is $\pi$-non-degenerate?\\
\end{que}

We expect that $\widetilde{d}_\alpha \neq \widetilde{d}^\alpha_{SH}$ in general (see Remark \ref{rmk:comparison to hofer metrics for ham} below).\\

\begin{ex}\label{ex:deltapm on ST*Tn} (strict, autonomous contact flows in $ST^*T^n$) (see also \cite[Theorem 5.8]{arl21}) Consider the unit cotangent bundle $ST^*T^n$ of the $n$-torus $T^n$. We use coordinates $(q,p) \in ST^*T^n \subseteq T^*T^n \approx T^n \times \R^n$ where $q \in T^n = ([0,1]/(0=1))^n$ and $p \in ST^*_qT^n \approx S^{n-1} \subseteq \R^n$. The contact form is the standard contact form given by $\alpha = pdq$. Let $H = H(p)$ and $G = G(p)$ be autonomous Hamiltonians depending only the $p$-coordinate. We claim that 	
\begin{equation}\label{eq:deltapm on ST^*T^N}
	\widetilde{\delta}^+_\alpha(f^H_t, f^G_t) = t \max_M (G-H), \quad \widetilde{\delta}^-_\alpha(f^H_t, f^G_t) = -t \min_M (G-H) \quad \text{for } t \geq 0.
\end{equation}
In particular, it follows that 
\begin{equation}
	\widetilde{d}^\alpha_{SH}(Id,f^H_t) = \widetilde{d}_\alpha(Id, f^H_t) = t \max_M |H|,
\end{equation}
which means that $f^H_t$ is length minimizing for $\widetilde{d}_{SH}$ and $\widetilde{d}_\alpha$ for all times.

To see this, we need to understand the relation $\lprec$ for contact isotopies of the above form. Let $H(p)$ and $G(p)$ be Hamiltonians as above. Then $f^H_t \lprec f^G_t$ if and only if $H(p) < G(p)$ for all $p \in S^{n-1}$. Indeed, first note that the contactomorphism $\phi^H_t$ is strict since $H$ is preserved by the Reeb flow $\phi^\alpha_t(q,p) = (q+tp,p)$. Thus, the contact isotopy $(\phi^H_t)^{-1} \circ \phi^G_t$ is generated by the Hamiltonian
\begin{equation}
F_t = -H \circ \phi^H_t + G \circ \phi^H_t,
\end{equation}
which is positive if $H < G$. Conversely, if $f^H_t \lprec f^G_t$, then in the notation of Section 3.2 in \cite{ep00}, we have that 
\begin{equation} 
H(p) = r_\pm(p,f^H_t) \leq r_\pm(p,f^G_t) = G(p)
\end{equation}
for any $p \in S^{n-1}$ (viewed now also as a non-zero cohomology class of $T^n$). In fact, if $f^H_t \lprec f^G_t$ for some fixed $t$, then also $f^H_t \lprec f^{G - \eps}_t$ for all sufficiently small $\eps > 0$, so we see that $H(p) \leq G(p) - \eps < G(p)$.

Since $s \mapsto \phi^H_{st}, s \in [0,1],$ and $s \mapsto \phi^\alpha_{s\eps} \phi^G_{st},s \in [0,1],$ are generated by the Hamiltonians $tH$ and $tG + \eps$, resp., it follows that $f^H_t \lprec f^\alpha_\eps f^G_t$ if and only if $t H < t G + \eps$, and $f^H_t \gprec f^\alpha_{-\eps} f^G_t$ if and only if $t H > tG - \eps$. In view of (\ref{eq:deltapm via cones}), this is equivalent to (\ref{eq:deltapm on ST^*T^N}) if $t \geq 0$.\\
\end{ex}

\begin{rmk}
One may ask about the continuity properties of the relation $\leq$ with respect to the interval topology $\widetilde{\mathcal{O}}_\lprec$, e.g. are the intervals
\begin{equation}
I^\leq_{f_0,f_1} \coloneqq \{f \in \tildecont(M,\xi)|\, f_0 \leq f \leq f_1\}
\end{equation}
or more generally the cones 
\begin{equation}
I^{\pm, \leq}_f \coloneqq \{ g \in \tildecont(M,\xi)|\, f^{\pm 1} \leq g^{\pm 1} \}
\end{equation}
closed with respect to $\widetilde{\mathcal{O}}_\lprec$? It is easy to see that $I^{\pm, \leq}_f \subseteq \overline{I^{\pm}_f}$ is always true, but it is unclear when equality holds.

Example \ref{ex:deltapm on ST*Tn} shows that for an autonomous Hamiltonian $H$ on $ST^*T^n$ that depends only on the $p$-coordinate, $f^H \in \overline{I^{\pm}_{id}}$ if and only if $f^H \in I^{\pm, \leq}_{id}$: If, say, $f^H \not\geq id$, then for some $p \in S^{n-1}$, $H(p) < 0$. It follows that for any $\eps > 0$ with $H(p) + \eps < 0$, $f^\alpha_\eps f \not\geq id$. In particular, $I^\lprec_{f^\alpha_{-\eps} f, f^\alpha_\eps f} \cap I^{+}_{id} = \emptyset$. Since the intervals $I^\lprec_{f^\alpha_{-\eps} f, f^\alpha_\eps f}$ form a neighbourhood basis of $f$ with respect to $\widetilde{\mathcal{O}}_\lprec$, this is equivalent to $f \not\in \overline{I^{+}_{id}}$. The case $f^H \not\leq id$ works similarly.

The same question can be asked about the relation $\leq$ and the interval topology on the other spaces considered below. In the case of the path-connected component of the zero section in the space of Legendrians in the $1$-jet bundle $J^1L$ of a closed manifold $L$, it follows in a similar way from Example \ref{ex:1-jets in J1L} below that for a function $f:L \to \R$, $j^1 f \in \overline{I^\pm_{j^1 0}}$ if and only if $f \in I^{\pm,\leq}_{j^1 0}$ if and only if $\pm f \geq 0$.\\
\end{rmk}

\subsection{Distances on $\Cont(M,\xi)$}\label{sec:cont(M,xi)}\hfill 

On $\Cont(M,\xi)$, one similarly defines $\delta^\pm_\alpha$ as before, except that we now take the infimum over \emph{all} contact isotopies with the given endpoints and not just those in a fixed homotopy class. The distances $\delta^\pm_\alpha$ then satisfy the analogue of Proposition \ref{prop:properties of deltapm} by the same proof, where in $(vi)$ we replace orderability by strong orderability. Thus, 
\begin{equation}
d_\alpha \coloneqq \max\{\delta^-_\alpha,\delta^+_\alpha,0\}
\end{equation}
defines a pseudo-metric on $\Cont(M,\xi)$ satisfying the analogues of $(i) - (v)$ in Proposition \ref{prop:properties of d_alpha}. By essentially the same proof as before, $\mathcal{O}_\lprec$ agrees with the metric topology of $d_\alpha$, and we have the following version of Theorem \ref{thm:non-deg of d_alpha} (see Section \ref{sec:proofs of main results} for the proof).

\begin{thm}\label{thm:non-degeneracy of d_alpha non-universal}
	On $\Cont(M,\xi)$, $d_\alpha$ either vanishes identically or is  non-degenerate. If $M$ is compact, then $d_\alpha$ vanishes identically if and only if there exists a positive loop of contactomorphisms on $M$. \eproof
\end{thm}

\begin{cor}
	Let $(M,\xi)$ be a closed contact manifold. Then the interval topology $\mathcal{O}_\lprec$ on $\Cont(M,\xi)$ is Hausdorff if and only if there does not exist a positive loop of contactomorphisms on $(M,\xi)$. \eproof\\
\end{cor}

Clearly, $d_\alpha$ is bounded by the Shelukhin-Hofer metric. Again, we do not know whether $d_\alpha$ is equivalent to the Shelukhin-Hofer metric on strongly orderable contact manifolds.\\

\begin{rmk}
	$\widetilde{d}_\alpha$ also induces the pseudo-metric
	\begin{equation}
	\widehat{d}_\alpha(\phi_0,\phi_1) \coloneqq \inf\{\widetilde{d}_\alpha(f_0,f_1)|\,\pi(f_0)=\phi_0, \pi(f_1) = \phi_1\} 
	\end{equation}
	on $\Cont(M,\xi)$ where $\pi:\tildecont(M,\xi) \to \Cont(M,\xi)$ denotes the covering map. It satisfies the analogues of properties $(i)-(v)$ in Proposition \ref{prop:properties of d_alpha}, and it is either non-degenerate or vanishes identically. But in contrast to $d_\alpha$, the proof of Theorem \ref{thm:non-deg of d_alpha} actually shows that $\widetilde{d}_\alpha$ vanishes identically if $\widehat{d}_\alpha$ is degenerate (see Remark \ref{rem:degeneracy of dhat implies vanishing of dtilde}). It follows that $\widehat{d}_\alpha$ is non-degenerate if and only if $(M,\xi)$ is orderable. In particular, $\widehat{d}_\alpha$ is different from $d_\alpha$ in general. For example, $ST^*S^2$ with its standard contact structure has periodic Reeb flow and is orderable \cite[Theorem 1.1]{cn16}.\\
	
	It follows directly from the definitions and Proposition \ref{prop:delta as integral} that
	\begin{equation}
	d_\alpha \leq \widehat{d}_\alpha \leq d^\alpha_{SH}.
	\end{equation}\\
	
	There is also a version of the interval topology which agrees with the metric topology of $\widehat{d}_\alpha$. Namely, it is generated by sets of the form
	\begin{equation}
	I^\lprec_{\phi_t} \coloneqq \bigcap\limits_{f} \pi(I^\lprec_{f,[\phi_t]f}),
	\end{equation}
	where $\phi_t \in \Cont(M,\xi), t \in [0,1],$ is a contact isotopy (not necessarily starting at the identity), and the intersection is taken over all $f \in \tildecont(M,\xi)$ with $\pi(f) = \phi_0$.\\
\end{rmk}

\begin{rmk}\label{rmk:comparison to hofer metrics for ham}
The definitions of $\delta^\pm_\alpha$, $d_\alpha$, and $\widehat{d}_\alpha$ are reminiscent of the semi-norms $\rho^\pm$, $\rho^+ + \rho^-$, and $\rho_f$  in Hofer's geometry (see \cite{ep93}, \cite{mcd02}). It was shown by McDuff \cite{mcd02} that $\rho_f$ and $\rho^+ + \rho^-$ may differ, so we expect that $\widehat{d}_\alpha \neq d_\alpha$ is possible even if $M$ is strongly orderable. \\
\end{rmk}

\begin{rmk}
Assume now that $(M,\xi)$ is compact and orderable, but not strongly orderable, i.e. $\widetilde{d}_\alpha$ and $\widehat{d}_\alpha$ define ($\pi$-)non-degenerate pseudo-metrics on $\tildecont(M,\xi)$ and $\Cont(M,\xi)$, respectively, and there exists a positive loop $\phi_t \in \Cont(M,\xi), t \in [0,1]$. By Lemma 3.1.A in \cite{ep00} we may assume that the generating Hamiltonian $F_t$ of $\phi_t$ is $1$-periodic. Denote by $f$ the element of $\tildecont(M,\xi)$ represented by the loop $\phi_t$. Following the constructions in \cite{fpr18}, we can define 
\begin{equation}
\begin{split}
\widetilde{\nu}^{f,+}_{FPR}(g) &\coloneqq \min\{ k \in \Z |\, f^k \geq g\},\\
\widetilde{\nu}^{f,-}_{FPR}(g) &\coloneqq \max\{ k \in \Z |\, f^k \leq g\},\\
\widetilde{\nu}^f_{FPR}(g) &\coloneqq \max \{|\widetilde{\nu}^{f,+}_{FPR}(g)|,|\widetilde{\nu}^{f,-}_{FPR}|(g)\}\\ 
&\,\,= \min\{k \in \Z_{\geq0}|\, f^{-k} \leq g \leq f^k\},
\end{split}
\end{equation}
for any $g \in \tildecont(M,\xi)$. If the Reeb flow of the contact form $\alpha$ is $1$-periodic and $f = f^\alpha_1$, this definition reduces to the norms $\nu_\pm$ and $\nu$ considered in \cite{fpr18}. By the same arguments that were used in \cite{fpr18} it is straightforward to check that $\widetilde{\nu}^f_{FPR}$ defines a pseudo-norm on $\tildecont(M,\xi)$ that is non-degenerate if and only if $(M,\xi)$ is orderable. Furthermore, $\widetilde{\nu}^{f,\pm}_{FPR}$ and $\widetilde{\nu}^f_{FPR}$ are conjugation invariant since $f$ is in the center of $\tildecont(M,\xi)$ and $\geq$ is bi-invariant. 

We claim that 
\begin{equation}\label{eq:comparison nu_FPR and d}
(\widetilde{\nu}^{f}_{FPR}(g) - 1) \widetilde{\tau}_\alpha([id_M],f)\leq \widetilde{d}_\alpha([id_M],g) \leq \widetilde{\nu}^{f}_{FPR}(g) \widetilde{d}_\alpha([id_M],f),
\end{equation}
where $\widetilde{\tau}_\alpha$ denotes the natural extension of $\tau_\alpha$ to the universal cover $\tildecont(M,\xi)$ (see Section \ref{sec:hedickes lorentzian distance function}).

Indeed, let $G_t$ be a Hamiltonian that generates $g$. If $\int\limits_0^1 \max G_t dt < k \int\limits_0^1 \min F_t dt$ for $k \geq 0$, then $f^k \gprec g$, and if $\int\limits_0^1 \min G_t dt > k \int\limits_0^1 \min F_t dt$ for $k \leq 0$, then $f^k \lprec g$. This implies the first inequality. On the other hand, if $f^k \geq g$ for $k \geq 0$, then it follows from Lemma 2.5 in \cite{fpr18} that $g$ can be generated by a Hamiltonian $G_t$ that satisfies $\int\limits_0^1 \max G_t dt \leq k \int\limits_0^1 \max F_t dt$, and if $f^k \leq g$ for $k \leq 0$, then $g$ can be generated by a Hamiltonian $G_t$ that satisfies $\int\limits_0^1 \min G_t dt \geq k \int\limits_0^1 \max F_t dt$. This implies the second inequality.

In the special case that the Reeb flow of the contact form $\alpha$ is $1$-periodic and $f = f^\alpha_1$, (\ref{eq:comparison nu_FPR and d}) simply reads 
\begin{equation}
\widetilde{\nu}^{f^\alpha_1}_{FPR}(g) - 1 \leq \widetilde{d}_\alpha([id_M],g) \leq \widetilde{\nu}^{f^\alpha_1}_{FPR}(g)
\end{equation}
since $\widetilde{\tau}_\alpha([id_M],f^\alpha_1) = \widetilde{d}_\alpha([id_M],f^\alpha_1) = 1$ (see Section \ref{sec:reeb flow}).

Note that the left-most and the right-most terms in (\ref{eq:comparison nu_FPR and d}) are invariant under conjugating $g$, whereas the middle term may change. In particular, these inequalities give bounds on $\widetilde{d}_\alpha([id_M],g)$ for all $g$ in a fixed conjugacy class.\\

Following \cite{fpr18}, we define
\begin{equation}
\nu^{f}_{FPR}(\psi) \coloneqq \min\{\widetilde{\nu}^{f}_{FPR}(g)|\, \pi(g) = \psi\}
\end{equation}
for $\psi \in \Cont(M,\xi)$, where the minimum is achieved since $\widetilde{\nu}^{f}_{FPR}(g) \geq 0$ is integer valued.

By taking the infimum over all $g$ with $\pi(g) = \psi$ in (\ref{eq:comparison nu_FPR and d}), we obtain
\begin{equation}\label{eq:comparison nu_FPR and dhat}
(\nu^{f}_{FPR}(\psi) - 1) \widetilde{\tau}_\alpha([id_M],f)\leq \widehat{d}_\alpha(id_M,\psi) \leq \nu^{f}_{FPR}(\psi) \widetilde{d}_\alpha([id_M],f).
\end{equation}

As before the right-most and the left-most terms are invariant under conjugating $\psi$, so that we get bounds on $\widehat{d}_\alpha(id_M,\psi)$ for all $\psi$ in a fixed conjugacy class.\\
\end{rmk}

\subsection{Distances on orbit spaces of subsets}\hfill

Let $A \subseteq M$ be a closed subset, and let
\begin{equation}
\mathscr{L}(A) \coloneqq \{\phi(A) \subseteq M|\, \phi \in \Cont(M,\xi) \} 
\end{equation}
be the orbit of $A$ under $\Cont(M)$. Denote by $i: A \hookrightarrow M$ the inclusion of $A$. An isotopy $i_t, t \in [0,1],$ starting at $i$ is a \emph{positive} (resp. \emph{non-negative}) isotopy if it is induced by an ambient contact isotopy $\phi_t$ with Hamiltonian $H_t$ that satisfies $H_t|_{\phi_t(A)} > 0$ (resp. $H_t|_{\phi_t(A)} \geq 0$). Note that the condition on $H_t$ does not depend on the choice of $\phi_t$ since $H_t \circ i_t = \alpha(\frac{d i_t}{dt})$, and that the existence of a positive or non-negative isotopy does not depend on the parametrization $i$ of $A$. As before, non-negative and positive isotopies induce  quasi-orders $\leq$ and (if $A$ is compact) $\lprec$ on $\mathscr{L}(A)$. $A$ is called \emph{strongly orderable} if $\leq$ defines a partial order. For compact $A_0, A_1 \in \mathscr{L}(A)$, we define the interval $I^\lprec_{A_0,A_1} \coloneqq \{B \in \mathscr{L}(A)|\, A_0 \lprec B \lprec A_1 \}$, and the collection of such intervals generates a topology $\mathcal{O}^A_\lprec$ on $\mathscr{L}(A)$.

The distances $\delta^{\pm,A}_\alpha$ induce distances on $\mathscr{L}(A)$ via
\begin{equation}
\begin{split}
\delta^{\pm,A}_\alpha(A_0,A_1) &\coloneqq \inf \{ \delta^\pm_\alpha(id_M, \phi)|\, \phi \in \Cont(M),~\phi(A_0) = A_1\},\quad A_0,A_1 \in \mathscr{L}(A),\\
d^A_\alpha(A_0,A_1) &\coloneqq \max \{\delta^{+,A}_\alpha(A_0,A_1),\delta^{-,A}_\alpha(A_0,A_1),0\}.
\end{split}
\end{equation}
It follows that $\delta^{\pm,A}_\alpha$ and $d^A_\alpha$ on $\mathscr{L}(A)$ satisfy the analogues of the properties in Propositions \ref{prop:properties of deltapm} and \ref{prop:properties of d_alpha}.\\

\begin{rmk}
One can also define the pseudo-metric 
\begin{equation}
\overline{d}^A_\alpha(A_0,A_1) \coloneqq \inf\{d_\alpha(id_M,\phi) |\, \phi(A_0)=A_1\}
\end{equation} 
induced by $d_\alpha$ on $\mathscr{L}(A)$. It is easy to check that it satisfies the properties in Proposition \ref{prop:properties of d_alpha} where we use the interval topology $\overline{\mathcal{O}}_\lprec$ generated by sets of the form
\begin{equation}
I_{\phi_0,\phi_1}^\lprec(A) \coloneqq \{\phi(A) |\, \phi_0 \lprec \phi \lprec \phi_1\}, \quad \phi_0,\phi_1 \in \Cont(M,\xi).
\end{equation}
An analogous remark also applies to the spaces $\mathscr{L}(i)$, $\widetilde{\mathscr{L}}(A)$ and $\widetilde{\mathscr{L}}(i)$ considered below.\\
\end{rmk}

Similar constructions also work for \emph{parametrized} subsets. Let 
\begin{equation}
	\mathscr{L}(i) \coloneqq \{\phi \circ i:A \hookrightarrow M|\, \phi \in \Cont(M,\xi)\}
\end{equation}
denote the orbit of $i$ under the action of $\Cont(M,\xi)$. The relations $\leq$ and (if $A$ is compact) $\lprec$ on $\mathscr{L}(i)$ are defined as before. If $A$ is compact, we define the intervals $I^\lprec_{i_0,i_1} \coloneqq \{j \in \mathscr{L}(i)|\, i_0 \lprec j \lprec i_1 \}$ for $i_0, i_1 \in \mathscr{L}(i)$, and the collection of such intervals generates a topology $\mathcal{O}^i_\lprec$ on $\mathscr{L}(i)$.

Then $\delta^{\pm,i}_\alpha$ and $d^i_\alpha$ on $\mathscr{L}(i)$ induced by $\delta^\pm_\alpha$ satisfy the analogues of the properties in Proposition \ref{prop:properties of d_alpha}.\\
  
In particular, this means that if $A$ is compact, the metric topologies of $d_\alpha^A$ and $d_\alpha^i$ agree with the interval topologies $\mathcal{O}^A_\lprec$ and $\mathcal{O}^i_\lprec$. The proof for $d_\alpha^i$ if $A$ is non-empty and  the proof for $d_\alpha^A$ if $A$ is non-empty and not equal to $M$ differs only notationally from the proof for $\Cont(M,\xi)$. If $A$ is empty then $\mathscr{L}(i)$ and $\mathscr{L}(A)$ both only consist of a single point, and if $A = M$, $\mathscr{L}(A)$ only consists of a single point, so the assertion is satisfied trivially.\\

Following \cite{rz18} (see also \cite{ush15}), we show in Section \ref{sec:proof of chekanov's dichotomoy for d_alpha} that $d^A_\alpha$ and $d^i_\alpha$ satisfy a version of Chekanov's dichotomy.

\begin{thm}\label{thm:chekanov's dichotomy for d_alpha}
	Let $i:L \hookrightarrow M^{2n+1}$ be a proper embedding of a submanifold $L \subseteq M$. If $d_\alpha^L$ is non-degenerate or if for all $i_0,i_1 \in \mathscr{L}(i)$ with $i_0(L) \neq i_1(L)$ we have that  $d_\alpha^i(i_0,i_1) > 0$, then $L$ is coisotropic. If $\dim L \leq n$, then $d_\alpha^L$ and $d_\alpha^i$ vanish identically unless $L$ has a Legendrian component. 
	
	If $L$ is a connected Legendrian submanifold, then either $d_\alpha^L$ vanishes identically or $d_\alpha^L$ is non-degenerate, and either $d_\alpha^i$ vanishes identically or for any $i_0,i_1 \in \mathscr{L}(i)$ with $i_0(L) \neq i_1(L)$, we have that $d_\alpha^i(i_0,i_1) > 0$. If in addition $L$ is compact, then $d_\alpha^L$ vanishes identically if and only if $L$ admits a positive loop, and $d_\alpha^i$ vanishes identically if and only if $i$ admits a positive loop.
\end{thm} 

As an immediate consequence, we get 

\begin{cor}
	Let $L \subseteq (M,\xi)$ be a closed, connected Legendrian. Then the interval topology $\mathcal{O}^L_\lprec$ on $\mathscr{L}(L)$ is Hausdorff if and only if $L$ does not admit a positive loop.\eproof\\
\end{cor}

\begin{rmk}
Note that if $i:L \hookrightarrow M$ is a proper Legendrian embedding, then $d_\alpha^i$ is never non-degenerate since isotopies of $i(L)$ are induced by ambient contact isotopies whose Hamiltonian vanishes along the image of $i(L)$. In this case, we can define an equivalence relation on $\mathscr{L}(i)$ where $i_0$ and $i_1$ are equivalent if and only if they are related by an isotopy of $L$ which starts at the identity. Denote the set of equivalence classes by $\widehat{\mathscr{L}}(i)$. Then it is easy to see that $d_\alpha^i$ descends to a pseudo-metric $\widehat{d}_\alpha^i$ on $\widehat{\mathscr{L}}(i)$ satisfying properties $(i)-(v)$ in Proposition \ref{prop:properties of d_alpha}. Similarly, the quasi-orders $\leq$ and (if $L$ is closed) $\lprec$ descend to quasi-orders on $\widehat{\mathscr{L}}(i)$, so we can again define an interval topology  $\widehat{\mathcal{O}}^i_\lprec$. As before, the metric topology induced by $\widehat{d}_\alpha^i$ is agrees with $\widehat{\mathcal{O}}^i_\lprec$ if $L$ is closed.
\end{rmk}

\begin{que}\label{que:dhat non-deg}
	Is $\widehat{d}_\alpha^i$ non-degenerate if it does not vanish identically?\\
\end{que}

Question \ref{que:dhat non-deg} seems to be hard to answer and its nature appears to be very different from the results discussed above. For example, assume that the answer is affirmative for the isotopy class of the zero-section in the 1-jet bundles $J^1L$ of a manifold $L$. Then this would imply that there does not exist a compactly supported Legendrian isotopy $i_t:L \to J^1L$ starting at the embedding of the zero-section so that $i_1(L) = i_0(L)$, when $i_1$ and $i_0$ represent different elements of the mapping class group of $L$, since the energy of such an isotopy can be made arbitrarily small by conjugation with the contactomorphism $(q,p,z) \mapsto (q,\lambda p, \lambda z)$ for small $\lambda > 0$. Conversely, if the answer is no, then there has to exist such an isotopy $i_t$.\\

It is clear that if $d_\alpha^i$ vanishes identically then also $d_\alpha^L$ vanishes identically as $d^i_\alpha(i_0,i_1) \geq d^L_\alpha(i_0(L),i_1(L))$.

\begin{que}
	If $d_\alpha^L$ vanishes identically, then does $d_\alpha^i$ vanish identically as well?\\
\end{que}

\begin{rmk}
Note that for the parametrized and unparametrized Shelukhin-Chekanov-Hofer metrics the answer to the corresponding question is yes. Indeed, assume that for two fixed Legendrians $L_0 \neq L_1$ and any $\eps > 0$, there exists $\phi_t \in \Cont(M,\xi), t \in [0,1],$ with $\phi_0 = id_M$, $\phi_1(L_0)= L_1$, and $\Vert \phi_t \Vert_\alpha < \eps$. Let $\psi_t \in \Cont(M,\xi), t \in [0,1], \psi_0 = id_M,$ be so that $\psi_1(L_0) \neq L_0$ and $\supp\, \psi_t \cap L_1 = \emptyset$. Then the concatenation $(\phi_t * (\psi_1 \circ \phi_{1-t}))|_{L_0}$ is well-defined and its energy is bounded by $(1 + \max \chi)\eps$, where $\chi$ is defined by $\psi_1^*\alpha = \chi \alpha$. Taking $\eps \to 0$ while fixing $\psi_t$, we obtain contact isotopies $\theta_t, t \in [0,1],$ starting at the identity which satisfy $\theta_1|_{L_0} = \psi_1|_{L_0}$ and have arbitrarily small energy. By the analogue of Theorem \ref{thm:chekanov's dichotomy for d_alpha} for the Shelukhin-Chekanov-Hofer metric (which in the unparametrized case is proven in \cite{rz18} and for parametrized Legendrians follows by the same proof that we use in Section \ref{sec:proof of chekanov's dichotomoy for d_alpha} for $d^i_\alpha$), it follows that the parametrized metric vanishes identically as well. 

This proof does however not work for $d^L_\alpha$ and $d^i_\alpha$ as we can only assume one-sided bounds on the generating Hamiltonian of $\phi_t$, so the Hamiltonian of $\psi_1 \circ \phi_{1-t}$ may not be small.\\
\end{rmk}

Assume that $i_L \hookrightarrow M$ is a proper Legendrian embedding of a connected manifold $L$. Like $d_\alpha$, $d_\alpha^L$ and $d_\alpha^i$ are bounded from above by the unparametrized or parametrized Shelukhin-Chekanov-Hofer metrics on $\mathscr{L}(L)$ and $\mathscr{L}(i)$, respectively. 

\begin{que}
	If $d_\alpha^L$ (resp. $d_\alpha^i$) does not vanish identically, then is it equivalent to the unparametrized (resp. parametrized) Shelukhin-Chekanov-Hofer metric?\\
\end{que}

\begin{ex}\label{ex:1-jets in J1L} (\emph{1-jets in $J^1L$})
	Consider the $1$-jet bundle $J^1 L$ of a closed manifold $L$ with the canonical contact form $\lambda = dz - p dq$. We identify $L$ with the zero section in $J^1 L$. Following Viterbo \cite{vit92}, one defines a function $c_-: \mathscr{L}(L) \to \R$, which is non-decreasing along non-negative isotopies by Lemma 5.2 in \cite{cn10} and on $1$-jets of functions $f:L \to \R$ it is equal to $c_-(j^1f) = \min f$ (see Example 3.1 in \cite{cn20}). Using this, we can compute $d^L_\alpha$ on $1$-jets. First note that if $i:L \to J^1L$ has a generating function $S$, then $\phi^\alpha_t \circ i$ has the generating function\footnote{Technically, $S+t$ is not quadratic at infinity, but we can first stabilize $S \to S + \xi^2, \xi \in \R,$ and then add a term of the form $\rho(\xi)t$, where $\rho: \R \to [0,1]$ is a compactly supported cut-off function which is equal to 1 on a large interval $[-C,C]$ and has sufficiently small derivatives. Then $S + \xi^2 + \rho(\xi)t$ is a generating function quadratic at infinity for $\phi^\alpha_t \circ i$ and $c_-(\phi^\alpha_t(i(L))) = c_-(i(L)) + t$.} $S + t$ where $\phi^\alpha_t$ denotes the time-$t$ Reeb flow on $J^1L$. From this it follows immediately that $c_-(\phi^\alpha_t(i(L))) = c_-(i(L)) + t$. Let $H_t: J^1L \to (-\eps,\infty)$ be a Hamiltonian, and denote its associated contact isotopy by $\phi^H_t$. Then $\phi^\alpha_{\eps t} \circ \phi^H_t$ is positive. Thus, $c_-(\phi^H_t(L')) = c_-(\phi^\alpha_{\eps t}(\phi^H_t(L'))) - \eps \geq c_-(L') - \eps$ for any $L' \in \mathscr{L}(L)$. This together with $c_-(j^1f) = \min f$ implies that $\delta^{-,L}_\alpha(L,j^1f) \geq - \min f$. Since the Hamiltonian $H(q,p,z) = f(q)$ generates the contact isotopy $(q,p,z) \mapsto (q,p+t df, z+tf)$, we actually have equality 
	\begin{equation}
	\delta^{-,L}_\alpha(L,j^1f) = - \min f.
	\end{equation} 
	
	Using invariance under left action by a strict contactomorphism, we obtain 
	\begin{equation}
	\delta^{+,L}_\alpha(L,j^1f) = \delta^{-,L}_\alpha(j^1f,L) = \delta^{-,L}_\alpha(L,j^1(-f)) = - \min(-f) = \max f.
	\end{equation} 
	It follows that $d_\alpha^L(L,j^1f) = \max |f| = d_{SCH}(L,j^1f)$. Again using right invariance this implies that $d_\alpha^L(j^1f,j^1g) = \max |f-g| = d_{SCH}(j^1f,j^1g)$.\\
\end{ex}

\begin{rmk}
We can also consider the interval topologies $\widetilde{\mathcal{O}}^i_\lprec$ and $\widetilde{\mathcal{O}}^A_\lprec$ and the pseudo-metrics $\widetilde{d}^i_\alpha$ and $\widetilde{d}^A_\alpha$ on the universal covers $\widetilde{\mathscr{L}}(i)$ and $\widetilde{\mathscr{L}}(A)$ of $\mathscr{L}(i)$ and $\mathscr{L}(A)$, respectively. By the same proofs as before (see Remark \ref{rmk:proof of chekanovs dichotomy on universal covers}) we have the following.

\begin{thm}
Let $i:L \hookrightarrow M^{2n+1}$ be the embedding of a properly embedded submanifold $L \subseteq M$. If $\widetilde{d}_\alpha^L$ on $\widetilde{\mathscr{L}}(L)$ is $\pi$-non-degenerate or if for all $\iota_0,\iota_1 \in \widetilde{\mathscr{L}}(i)$ with $\pi(\iota_0)(L) \neq \pi(\iota_1)(L)$ we have that $\widetilde{d}_\alpha^i(\iota_0,\iota_1) > 0$, then $L$ is coisotropic. If $\dim L \leq n$, then $\widetilde{d}_\alpha^L$ and $\widetilde{d}_\alpha^i$ vanish identically unless $L$ has a Legendrian component. 
	
	If $L$ is a connected Legendrian submanifold, then either $\widetilde{d}_\alpha^L$ vanishes identically or $\widetilde{d}_\alpha^L$ is $\pi$-non-degenerate, and either $\widetilde{d}_\alpha^i$ vanishes identically or for any $\iota_0,\iota_1 \in \widetilde{\mathscr{L}}(i)$ with $\pi(\iota_0)(L) \neq \pi(\iota_1)(L)$, $\widetilde{d}_\alpha^i(\iota_0,\iota_1) > 0$. If in addition $L$ is compact, then $\widetilde{d}_\alpha^L$ vanishes identically if and only if $L$ admits a contractible positive loop, and $\widetilde{d}_\alpha^i$ vanishes identically if and only if i admits a contractible positive loop. \eproof
\end{thm}

\begin{cor}
	Let $L \subseteq (M,\xi)$ be a closed, connected Legendrian. Then the interval topology $\mathcal{O}^L_\lprec$ on $\widetilde{\mathscr{L}}(L)$ is $\pi$-Hausdorff if and only if $L$ does not admit a contractible positive loop. \eproof\\
\end{cor}
\end{rmk}

\subsection{Hedicke's Lorentzian distance function}\label{sec:hedickes lorentzian distance function}\hfill

The definitions in \cite{hed21} of the Lorentzian distance functions $\tau_\alpha$ and $\tau^A_\alpha$ on $\Cont(M,\xi)$ and $\mathcal{L}(A)$ for a closed Legendrian $A$ extend naturally also to $\mathcal{L}(A)$ and $\mathcal{L}(i)$ for inclusions $i:A \to M$ of arbitrary closed subsets and to the universal covers $\tildecont(M,\xi)$, $\widetilde{\mathcal{L}}(A)$, and $\widetilde{\mathcal{L}}(i)$ of these spaces. These will be denoted by $\tau^A_\alpha$, $\tau^i_\alpha$, $\widetilde{\tau}_\alpha$, $\widetilde{\tau}^A_\alpha$, and $\widetilde{\tau}^i_\alpha$.\\

These different versions of the Lorentzian distance functions can be expressed in terms of the different versions of $\delta^+_\alpha$. For this, recall the following lemma from \cite{hed21}\footnote{The lemma in \cite{hed21} was stated in a weaker form not fixing the homotopy class, but it is clear from its proof that one may also assume that $\phi^G_t$ and $\phi^H_t, t \in [0,1],$ are homotopic relative to the endpoints.}.

\begin{lem}\label{lem:realize int min}
	Let $H_t:M \to \R,t \in [0,1],$ be a Hamiltonian on a closed contact manifold $(M,\ker \alpha)$. Then for any $\eps > 0$ there exists another Hamiltonian $G_t$ such that $[\phi^G_t] = [\phi^H_t] \in \tildecont(M,\xi)$ and 
	\begin{equation}
	\min\limits_{M} G_t \in \left(\int_0^1 \min_M H_s ds - \eps, \int_0^1 \min_M H_s ds + \eps \right).
	\end{equation} \\
\end{lem}

With this lemma, we can rewrite the definitions of $\widetilde{\tau}_\alpha$ and $\widetilde{\delta}^\pm_\alpha$ if $M$ is compact:

\begin{equation}\label{eq:tau,delta in terms of intergrals}
\begin{split}
\widetilde{\delta}^-_\alpha(f_0,f_1) &= \inf \left\{-\int_0^1 \min_M H_t dt |\, f^H f_0 = f_1 \right\} \\
&= -\sup \left\{\int_0^1 \min_M H_t dt |\, f^H f_0 = f_1 \right\}\\
\widetilde{\delta}^+_\alpha(f_0,f_1) &= \inf \left\{\int_0^1 \max_M H_t dt |\, f^H f_0 = f_1 \right\}, \\
\widetilde{\tau}_\alpha(f_0,f_1) &= \max \left\{ \sup \left\{\int_0^1 \min_M H_t dt |\, f^H f_0 = f_1 \right\},0 \right\}.
\end{split}
\end{equation}

In particular, we see that
\begin{equation}\label{eq:tau in terms of delta}
\widetilde{\tau}_\alpha = \max \{-\widetilde{\delta}^-_\alpha,0\},
\end{equation}
using that $\widetilde{\delta}^-_\alpha(f_0,f_1) < 0$ if and only if $f_0 \lprec f_1$. The same formulas hold also for $\tau_\alpha$ and $\delta^\pm_\alpha$. For $\tau^A_\alpha$, $\tau^i_\alpha$ $\delta^{\pm,A}_\alpha$, $\delta^{\pm,i}_\alpha$, $\widetilde{\tau}^A_\alpha$, $\widetilde{\tau}^i_\alpha$, $\widetilde{\delta}^{\pm,A}_\alpha$, and $\widetilde{\delta}^{\pm,i}_\alpha$, the corresponding formulas where the $\min$/$\max$ is taken over $A$ instead hold under the weaker assumption that the Reeb flow in $(M,\alpha)$ is complete, as is clear from the proof of Lemma 3.1 in \cite{hed21}.\\

Using this observation, it follows that Theorems 2.7 and 5.8 in \cite{hed21} hold for the constant $C_\alpha = 1$ whenever $(M,\xi)$ is closed and orderable, resp. when the Reeb flow is complete and the Legendrian is orderable.

\begin{prop}\label{prop:tau 1-lipschitz}
	If $\widetilde{d}_\alpha$ is non-degenerate, then $\widetilde{\delta}^\pm_\alpha$ is 1-Lipschitz with respect to $\widetilde{d}_\alpha$ in both arguments. If in addition $M$ is closed\footnote{If we apply the definition of $\widetilde{\tau}_\alpha$ to open manifolds, then it necessarily vanishes identically. In particular, it is 1-Lipschitz.}, then $\widetilde{\tau}_\alpha$ is 1-Lipschitz with respect to $\widetilde{d}_\alpha$ in both arguments. In particular, $\widetilde{\tau}_\alpha(f_0,f_1) \leq \widetilde{d}_\alpha(f_0,f_1) \leq \widetilde{d}_{SH}(f_0,f_1)$ for any $f_0,f_1 \in \tildecont(M,\xi)$.
\end{prop}

\bproof By (\ref{eq:tau in terms of delta}) and property (ii) in Proposition \ref{prop:properties of deltapm}, it is enough to prove the claim for $\widetilde{\delta}^-_\alpha$. We only show that $\widetilde{\delta^-_\alpha}$ is 1-Lipschitz in the second argument, the proof for the first argument is analogous. \\

For this, let $f, g, \widetilde{g} \in \tildecont(M,\xi)$ be so that $\widetilde{d}_\alpha(g,\widetilde{g}) = \eps > 0$. Denote $a \coloneqq \widetilde{\delta}^-_\alpha(f,\widetilde{g})$ and $b \coloneqq \widetilde{\delta}^-_\alpha(f,g)$. We will show that $|a-b| \leq \eps$. Without loss of generality, we may assume that $a \leq b$. Then we see that
\begin{equation}
b = \widetilde{\delta}^-_\alpha(f,g) \leq \widetilde{\delta}^-_\alpha(f,\widetilde{g}) + \widetilde{\delta}^-_\alpha(\widetilde{g},g) \leq \widetilde{\delta}^-_\alpha(f,\widetilde{g}) + \widetilde{d}_\alpha(\widetilde{g},g) = a + \eps,
\end{equation}
or equivalently, $|b-a| = b-a \leq \eps$.

The second part of the proposition is an immediate consequence of the first part since $\widetilde{\tau}(f_0,f_0) = 0$ by the orderability of $(M,\xi)$.
\eproof\\

The corresponding result applies with the same proof to the other variants of $\tau_\alpha$ and $\delta^\pm_\alpha$.\\

\begin{rmk}
The above results imply that for a closed, strongly orderable Legendrian $L$, $(\mathscr{L}(L),d^L_\alpha,\lprec,\leq,\tau^L_\alpha)$ is a strongly causal Lorentzian pre-length space in the sense of \cite{ks18}. However, it is in general not a Lorentzian length space. Indeed, any Lorentzian length space is by definition localizable (see \cite{ks18} for the definitions). In particular, if $(\mathscr{L}(L),d^L_\alpha,\lprec,\leq,\tau^L_\alpha)$ were a Lorentzian length space, every $L \in \mathscr{L}(L)$ would have an open neighbourhood $U$ (with respect to $d^L_\alpha$) together with a number $C>0$ so that for any non-negative isotopy $L_t \subseteq U$, we have that $L^{d^L_\alpha}(L_t) < C$, where $L^{d^L_\alpha}(L_t)$ denotes the length of $L_t$ viewed as a rectifiable curve in the metric space $(\mathscr{L}(L),d^L_\alpha)$. This is not satisfied already for the zero section in the 1-jet space of any closed manifold $L$. In fact, even any $C^k$-small neighbourhood of the zero section contains arbitrarily long non-negative isotopies. To see this, choose an open cover $\{U_i\}_{i = 1}^{N}, N \in \N,$ of $L$ so that for all $j \in \{1,...,N\}$, $U_j \setminus \bigcup_{i \neq j} U_i$ is non-empty. Pick a partition of unity $\{\lambda_i\}_{i = 1}^{N}$ subordinate to this cover. Note that $\lambda_i^{-1}(\{1\})$ is non-empty. Then for any $\eps > 0$ and $k \in \N_0$, the non-negative Legendrian isotopies given by the $1$-jets
\begin{equation}
j^1 \left(\eps \left(t \lambda_j + \sum_{i=1}^{j-1} \lambda_i + k\right) \right), \quad j \in \{1,...,N\},\, t \in [0,1],
\end{equation}
have length $\eps$ by Example \ref{ex:1-jets in J1L} above. Fix some $\delta > 0$. Let $\eps > 0$ be so small that for any $t \in [0,1]$, $j^1(t \eps)$ is $\frac{\delta}{2}$-close to the zero section in the $C^k$-norm. Then for any sufficiently large number $K \in \N$, the Legendrian isotopies
\begin{equation}
j^1 \left(\frac{\eps}{K} \left(t \lambda_j + \sum_{i=1}^{j-1} \lambda_i + k \right) \right), \quad j \in \{1,...,N\},\, k \in \{0,...,K-1\},\, t \in [0,1],
\end{equation}
are $\delta$-close to the zero section in the $C^k$-norm. Then the Legendrian isotopy $L_t, t \in [0,NK],$ from $j^1 0$ to $j^1 \eps$ which is equal to 
\begin{equation}
j^1 \left(\frac{\eps}{K} \left((t - (j-1+kN)) \lambda_j + \sum_{i=1}^{j-1} \lambda_i + k \right) \right) 
\end{equation}
for $t \in [j-1 + kN, j + kN], j \in \{1,...,N\}, k \in \{0,...,K-1\},$ is $\delta$-close to the zero section in the $C^k$-norm and has length $\frac{\eps}{K}NK = \eps N$. Taking the limit $N \to \infty$ (while potentially increasing $K$ as well), we find arbitrarily large non-negative isotopies in any $C^k$-neighbourhood of the zero section as claimed.

As any Legendrian admits a Weinstein neighbourhood, we expect that $(\mathscr{L}(L),d^L_\alpha,\lprec,\leq,\tau^L_\alpha)$ is not localizable for any closed, strongly orderable Legendrian $L$ in any contact manifold.\\

Similarly, for a closed, strongly orderable contact manifold $(M,\xi)$ with contact form $\alpha$, $(\Cont(M,\xi),d_\alpha,\lprec,\leq,\tau_\alpha)$ is a strongly causal Lorentzian pre-length space. Again, we expect it to be non-localizable. One reason is that the graph of a contactomorphism in $M \times M \times \R$ is Legendrian, and any Legendrian $C^1$-close to a graph is again graphical, so the above construction performed in a Weinstein neighbourhood of the graph gives a contact isotopy of arbitrarily large energy $C^k$-close to the identity.\\
\end{rmk}

\subsection{Reeb flow}\label{sec:reeb flow}\hfill

Generalizing Corollaries 2.4 and 2.5 in \cite{hed21}, we observe that the Reeb flow is always a ``length-minimizing geodesic" for $\widetilde{\delta}^\pm_\alpha$ if $(M,\xi)$ is orderable. In fact, as an immediate consequence of (\ref{eq:deltapm via cones}), the corresponding formula on $\Cont(M,\xi)$, and the left-invariance of $\lprec$, we have the following result.

\begin{prop}\label{prop:deltapm and reeb flow}
	Let $f \in \tildecont(M,\xi)$, $\phi \in \Cont(M,\xi)$, and assume that $M$ is closed. If $(M,\xi)$ is orderable, then 
	\begin{equation}
	\widetilde{\delta}^\pm_\alpha(Id_M,f^\alpha_t f) = \widetilde{\delta}^\pm_\alpha(Id_M,f) \pm t.
	\end{equation}
	
	If $(M,\xi)$ is strongly orderable, 
	\begin{equation}
	\delta^\pm_\alpha(Id_M,\phi^\alpha_t \phi) = \delta^\pm_\alpha(Id_M,\phi) \pm t. 
	\end{equation} \eproof\\
\end{prop}



Similarly, the following holds.

\begin{prop}
	Assume that the Reeb flow on $(M,\alpha)$ is complete. Let $i:A \to M$ be the embedding of a closed subset, and let $j \in \mathcal{L}(i)$, $\widetilde{j} \in \widetilde{\mathcal{L}}(i)$, $B \in \mathcal{L}(A)$, and $\widetilde{B} \in \widetilde{\mathcal{L}}(A)$. If $i$ is orderable, then
	\begin{equation}
	\widetilde{\delta}^{\pm,i}([i], f^\alpha_t \widetilde{j}) = \widetilde{\delta}^{\pm,i}([i],\widetilde{j}) \pm t.
	\end{equation}
	If $i$ is strongly orderable, then
	\begin{equation}
	\delta^{\pm,i}(i, \phi^\alpha_t \circ j) = \delta^{\pm,i}(i,j) \pm t.
	\end{equation}
	If $A$ is orderable, then
	\begin{equation}
	\widetilde{\delta}^{\pm,A}([A], f^\alpha_t \widetilde{B}) = \widetilde{\delta}^{\pm,A}([A],\widetilde{B})\pm t.
	\end{equation}
	If $A$ is strongly orderable, then
	\begin{equation}
	\delta^{\pm,A}(A, \phi^\alpha_t(B)) = \delta^{\pm,A}(A, B) \pm t.
	\end{equation} \eproof\\
\end{prop}

\subsection{The relative growth of contactomorphisms} \hfill

The relative growth of two contact isotopies can be estimated in terms of generating Hamiltonians of the two isotopies. 

\begin{prop}\label{prop:deltapm and relative growth}
	Assume that $M$ is compact and orderable, and let $f,g \in \tildecont(M,\xi)$ with $f \gprec Id_M$. If $\widetilde{\delta}^+_\alpha(Id_M,g) \geq 0$, then
\begin{equation}
\gamma(f,g) \leq \frac{\widetilde{\delta}^+_\alpha(Id_M,g)}{\widetilde{\tau}_\alpha(Id_M,f)},
\end{equation}

if $\widetilde{\delta}^+_\alpha(Id_M,g) < 0$, then
\begin{equation}
\gamma(f,g) \leq \frac{\widetilde{\delta}^+_\alpha(Id_M,g)}{\widetilde{\delta}^+_\alpha(Id_M,f)},
\end{equation}

if $\widetilde{\delta}^-_\alpha(Id_M,g) < 0$, then
\begin{equation}
\gamma(f,g) \geq \frac{\widetilde{\tau}_\alpha(Id_M,g)}{\widetilde{\delta}^+_\alpha(Id_M,f)},
\end{equation}

and if $\widetilde{\delta}^-_\alpha(Id_M,g) \geq 0$, then
\begin{equation}
\gamma(f,g) \geq  \frac{- \widetilde{\delta}^-_\alpha(Id_M,g)}{\widetilde{\tau}_\alpha(Id_M,f)}. 
\end{equation}
\end{prop}

\bproof Let $F_t$ and $G_t$ denote generating Hamiltonians of $f$ and $g$ so that $\int_0^1 \min\limits_M F_t dt > 0$. Assume that for $k \in \N$ and $p \in \Z_{\geq 0}$, we have that 
\begin{equation}
\frac{p}{k} > \frac{\int_0^1 \max\limits_M G_t dt}{\int_0^1 \min\limits_M F_t dt}.
\end{equation}
Then $f^{-p} g^k$ can be generated by a Hamiltonian $H_t$ which satisfies 
\begin{equation}
\int_0^1 \max\limits_M H_t dt = k  \int_0^1 \max\limits_M G_t dt - p \int_0^1 \min\limits_M F_t dt < 0,
\end{equation}
i.e. $Id_M \gprec f^{-p} g^k$, or $f^p \gprec g^k$. This shows that 
\begin{equation}
\gamma(f,g) \leq \max \left\{0,\inf\limits_{F,G} \frac{\int_0^1 \max\limits_M G_t dt}{\int_0^1 \min\limits_M F_t dt} \right\},
\end{equation}
which in the case that $\widetilde{\delta}^+_\alpha(Id_M,g) \geq 0$, means that 
\begin{equation}
\gamma(f,g) \leq \frac{\widetilde{\delta}^+_\alpha(Id_M,g)}{\widetilde{\tau}_\alpha(Id_M,f)}.
\end{equation}\\

Now assume that $\widetilde{\delta}^+_\alpha(Id_M,g) < 0$, or equivalently $g \lprec Id_M$, and that for $k \in \N$ and $p \in \Z_{<0}$,
\begin{equation}
\frac{p}{k} > \frac{\int_0^1 \max\limits_M G_t dt}{\int_0^1 \max\limits_M F_t dt}.
\end{equation}
Then $f^{-p} g^k$ can be generated by a Hamiltonian $H_t$ which satisfies 
\begin{equation}
\int_0^1 \max\limits_M H_t dt = k  \int_0^1 \max\limits_M G_t dt - p \int_0^1 \max\limits_M F_t dt < 0,
\end{equation}
which implies that $f^p \gprec g^k$. It follows that 
\begin{equation}
\gamma(f,g) \leq \inf\limits_{F,G} \frac{\int_0^1 \max\limits_M G_t dt}{\int_0^1 \max\limits_M F_t dt},
\end{equation}
which implies that 
\begin{equation}
\gamma(f,g) \leq \frac{\widetilde{\delta}^+_\alpha(Id_M,g)}{\widetilde{\delta}^+_\alpha(Id_M,f)}.
\end{equation}\\

On the other hand assume that for some $p \in \Z_{\geq 0}$ and $k \in \N$, we have that
\begin{equation}
\frac{p}{k} < \frac{\int_0^1 \min\limits_M G_t dt }{\int_0^1 \max\limits_M F_t dt}.
\end{equation}
Then $g^{-k} f^p$ can be generated by a Hamiltonian $H_t$ that satisfies
\begin{equation}
\int_0^1 \max\limits_M H_t dt = p  \int_0^1 \max\limits_M F_t dt - k \int_0^1 \min\limits_M G_t dt < 0,
\end{equation}
which implies that $f^p \lprec g^k$. This shows that 
\begin{equation}
\gamma(f,g) \geq \sup\limits_{F,G} \frac{\int_0^1 \min\limits_M G_t dt }{\int_0^1 \max\limits_M F_t dt}
\end{equation}
whenever $\int_0^1 \min\limits_M G_t dt > 0$, which means that 
\begin{equation}
\gamma(f,g) \geq \frac{\widetilde{\tau}_\alpha(Id_M,g)}{\widetilde{\delta}^+_\alpha(Id_M,f)}.
\end{equation}\\

If
\begin{equation}
\frac{p}{k} < \frac{\int_0^1 \min\limits_M G_t dt }{\int_0^1 \min\limits_M F_t dt}.
\end{equation}
for some $p \in \Z_{<0}$ and $k \in \N$, then $g^{-k} f^p$ can be generated by a Hamiltonian $H_t$ that satisfies
\begin{equation}
\int_0^1 \max\limits_M H_t dt = p  \int_0^1 \min\limits_M F_t dt - k \int_0^1 \min\limits_M G_t dt < 0,
\end{equation}
which again implies that $f^p \lprec g^k$. It follows that 
\begin{equation}
\gamma(f,g) \geq \min \left\{0,\sup\limits_{F,G} \frac{\int_0^1 \min\limits_M G_t dt }{\int_0^1 \min\limits_M F_t dt} \right\},
\end{equation}
where the supremum is taken over all generating Hamiltonians so that $\int_0^1 \min\limits_M F_t dt > 0$. In the case that $\widetilde{\delta}^-_\alpha(Id_M,g) \geq 0$ this implies that 
\begin{equation}
\gamma(f,g) \geq  \frac{- \widetilde{\delta}^-_\alpha(Id_M,g)}{\widetilde{\tau}_\alpha(Id_M,f)}. 
\end{equation}
\eproof\\

For autonomous Hamiltonians $F,G$ on $ST^*T^n$ which depend only on the coordinate on the fibre (see Example \ref{ex:deltapm on ST*Tn} above) Eliashberg and Polterovich \cite{ep00} showed that 
\begin{equation}
\gamma(f^F_t,f^G_t) = \max\limits_{p \in S^{n-1}} \frac{G(p)}{F(p)},
\end{equation}
whenever $F(p) > 0$ for all $p$ and $G(p) > 0$ for some $p$. In fact, this remains true even if $G(p) \leq 0$ for all $p$. Indeed, $(f^G_t)^k$ is generated by $k G$ and $(f^F_t)^q$ is generated by $q F$, so by the characterization of $\lprec$ given in Example \ref{ex:deltapm on ST*Tn}, $(f^F_t)^q \gprec (f^G_t)^k$ if and only if $q F(p) > k G(p)$ for all $p \in S^{n-1}$.

In this case the above inequalities, say in the case that $G > 0$ (i.e. $\widetilde{\delta}^+(Id_M,g) > 0$ and  $\widetilde{\delta}^-(Id_M,g) < 0$), read
\begin{equation}
\frac{\min G}{\max F} \leq \max \frac{G}{F} \leq \frac{\max G}{\min F},
\end{equation}
so we see that the inequalities cannot be replaced by equalities in general.

If on the other hand $f = f^\alpha_t, t > 0,$ and $g = f^\alpha_s, s \in \R$, then in view of Proposition \ref{prop:deltapm and reeb flow} the inequalities become 
\begin{equation}
\frac{s}{t} \leq \gamma(f^\alpha_t,f^\alpha_s) \leq \frac{s}{t},
\end{equation}
which means that they are, in fact, equalities.\\

\section{Reparametrizations of Hamiltonians}\label{sec:some proofs}

In this section we show the following lemma, which allows us to use reparametrizations of contact isotopies to get bounds on the maximum or minimum of the generating Hamiltonian.

\begin{lem}\label{lem:reparam of hamiltonian to make max const}
	Let $\lambda:[0,1] \to [0,\infty)$ be a continuous function, and let $\eps > 0$. Then there exists a diffeomorphism $\sigma: [0,1] \to [0,1]$ with $\sigma(0) = 0$ so that 
	\begin{equation}
	\int\limits_0^1 \lambda(t) dt + \eps >  \underset{t}{\max} \, \left(\lambda(\sigma(t))\sigma'(t) \right).
	\end{equation} 
\end{lem}

\bproof Denote $l \coloneqq \int\limits_0^1 \lambda(t) dt$. For $\delta > 0$, let $U \coloneqq \{t \in [0,1]| \frac{\lambda(t)}{l} < \delta \}$. Let $\phi: [0, \infty) \to [0, \infty)$ be a smooth function so that $\phi(x) = x$ for $x \leq \frac{1}{\delta}$ and $\phi(x) \in [\frac{1}{\delta},\frac{1}{\delta} + 1]$ for $x > \frac{1}{\delta}$.

There exists a number $\tau > 0$ and a diffeomorphism $\widetilde{\sigma}: [0,\tau] \to [0,1], \tau > 0,$ which are uniquely determined by the differential equation $\widetilde{\sigma}'(t) = \phi(\frac{l}{\lambda(\widetilde{\sigma}(t))}) \geq \min\{ \frac{1}{\delta}, \frac{l}{\max \lambda}\} > 0$ with the initial condition $\widetilde{\sigma}(0) = 0$. 

Note that 
\begin{equation}
\int\limits_0^1 \frac{1}{\phi(\frac{l}{\lambda(t)})} dt = \int\limits_0^{\tau}  \frac{1}{\phi(\frac{l}{\lambda(\widetilde{\sigma}(t))})} \widetilde{\sigma}'(t) dt = \tau.
\end{equation}

Note further that 
\begin{equation}
\begin{split}
\tau &= \int\limits_0^1 \frac{1}{\phi(\frac{l}{\lambda(t)})} dt = \int\limits_{[0,1] \setminus U} \frac{\lambda(t)}{l} dt + \int\limits_U \frac{1}{\phi(\frac{l}{\lambda(t)})} dt \\&= 1 + \int\limits_U \frac{1}{\phi(\frac{l}{\lambda(t)})} - \frac{\lambda(t)}{l}dt \in \left[1 - \delta \int\limits_U1 dt, 1 + \delta \int\limits_U1 dt \right]
\end{split}
\end{equation}
by definition of $U$ and $\phi$. In particular, $\tau \rightarrow 1$ as $\delta \rightarrow 0$. Now we define the diffeomorphism $\sigma: [0,1] \to [0,1], \sigma(t) = \widetilde{\sigma}(\frac{t}{\tau})$. Then it follows that $\sigma'(t) = \frac{1}{\tau} \phi(\frac{l}{\lambda(\sigma(t))})$, which implies that $\lambda(\sigma(t)) \sigma'(t) = \lambda(\sigma(t)) \frac{1}{\tau} \phi(\frac{l}{\lambda(\sigma(t))})$. If $\sigma(t) \in [0,1] \setminus U$, this is equal to  $\frac{l}{\tau}$, and if $\sigma(t) \in U$, $\lambda(\sigma(t)) \in [0,\delta l)$ and  $\phi(\frac{l}{\lambda(\sigma(t))}) \in [\frac{1}{\delta},\frac{1}{\delta} + 1]$ which implies that $\lambda(\sigma(t)) \sigma'(t) \in [0,(1+\delta)\frac{l}{\tau})$. In particular, for any $\eps > 0$, $\lambda(\sigma(t)) \sigma'(t) < l + \eps$ for all $t \in [0,1]$ if $\delta$ is sufficiently small. \eproof \\

\section{Proofs of Theorems \ref{thm:non-deg of d_alpha} and \ref{thm:non-degeneracy of d_alpha non-universal}}\label{sec:proofs of main results}

\subsection{Proof of the first part of Theorem \ref{thm:non-deg of d_alpha}} \hfill

To prove the first part of Theorem \ref{thm:non-deg of d_alpha}, assume that $\widetilde{d}_\alpha$ is not $\pi$-non-degenerate, i.e. by right-invariance of $\widetilde{d}_\alpha$, there exists $h \in \tildecont(M,\xi), \pi(h) \neq id_M,$ with $\widetilde{d}_\alpha(h, id_M) = 0$. Let $U \subseteq M$ be a open subset so that $\pi(h)(U) \cap U = \emptyset$. Using an argument by Eliashberg and Polterovich \cite{ep93} (see also \cite[Remark 7]{she17}), it is straightforward to show that for any $f \in \tildecont(M,\xi)$ with support in $U$, $\widetilde{d}_\alpha(f, id_M) = 0$. By the naturality of $\widetilde{d}_\alpha$ under conjugation, there exists for any point $x \in M$ a neighbourhood $U_x$ and $h_x \in \tildecont(M,\xi)$ so that $\widetilde{d}_\alpha(h_x,id_M)=0$ and $\pi(h_x)(U_x) \cap U_x = \emptyset$. By the fragmentation lemma for contact isotopies (see \cite{ban97}) and right-invariance, it then follows that $\widetilde{d}_\alpha$ vanishes identically. 

To prove the claim, let $h$ and $U$ be as above. For $1 \leq r \leq \infty$, we denote by $\textnormal{Cont}_0^r(U,\xi)$ the identity component of the group of $C^r$-contactomorphisms with compact support in $U$ and by $\tildecont^r(U,\xi)$ its universal cover (see \cite{tsu08} for the definitions). Observe that the definition of the pseudo-metrics extends naturally to $\tildecont^r(U,\xi)$. We first show that $\widetilde{d}_\alpha$ vanishes on commutators of elements in $\tildecont^r(U,\xi)$ for any $r \geq 1$. Let $f_1, f_2 \in \tildecont^r(U,\xi)$. We define $g \coloneqq [h, f_1^{-1}] \coloneqq f_1 h^{-1} f_1^{-1} h$ as the commutator of $h$ and $f_1^{-1}$. Then an easy argument shows that $[g,f_2] = [f_1, f_2]$. Furthermore,
\begin{equation}
\begin{split}
\widetilde{d}_\alpha([g, f_2],id_M) &= \widetilde{d}_\alpha(f_2^{-1} g^{-1} f_2 g, id_M) \leq \widetilde{d}_\alpha(f_2^{-1} g^{-1} f_2, id_M) + \widetilde{d}_\alpha(g, id_M) \\
&\leq (1 + C(\pi(f_2^{-1}))) \widetilde{d}_\alpha(g, id_M) = (1 + C(\pi(f_2^{-1}))) \widetilde{d}_\alpha(f_1 h^{-1} f_1^{-1} h, id_M)\\
&\leq (1 + C(\pi(f_2^{-1})))(1 + C(\pi(f_1))) \widetilde{d}_\alpha(h, id_M) = 0,
\end{split}
\end{equation}
where for any contactomorphism $\psi \in \textnormal{Cont}_0^r(M,\alpha)$, $C(\psi) \coloneqq \max \chi$, where $\chi:M \to (0, \infty)$ denotes the conformal factor of $\psi$ defined via $\psi^* \alpha = \chi \alpha$. In particular, $\widetilde{d}_\alpha([f_1, f_2],id_M) = 0$. Since $\tildecont^1(U,\xi)$ is perfect \cite{tsu08}, it follows that $\widetilde{d}_\alpha(-,id_M)$ vanishes on $\tildecont^1(U,\xi)$. Since paths in $\tildecont^1(U,\xi)$ can be $C^1$-approximated by paths in $\tildecont(U,\xi)$ and the generating Hamiltonian depends continuously (in the $C^1$-topology) on the path of contactomorphisms, $\widetilde{d}_\alpha$ also vanishes on the image of $\tildecont(U,\xi)$ in $\tildecont(M,\xi)$, as was claimed. \eproof\\

\begin{rmk}\label{rem:degeneracy of dhat implies vanishing of dtilde}
Above we proved that for any $f_1,f_2,h \in \tildecont^r(M,\xi)$ so that 
\begin{equation}
(\supp(f_1) \cup \supp(f_2)) \cap \pi(h)(\supp(f_1) \cup \supp(f_2)) = \emptyset,
\end{equation}
we have that 
\begin{equation}
\widetilde{d}_\alpha([f_1, f_2],id_M) \leq (1 + C(\pi(f_2^{-1})))(1 + C(\pi(f_1))) \widetilde{d}_\alpha(h, id_M).
\end{equation}
Let $U \subseteq M$ be so that $\pi(h)(U) \cap U = \emptyset$. Taking the infimum over all $h$ with fixed $\phi \coloneqq \pi(h)$, we obtain that 
\begin{equation}
\widetilde{d}_\alpha([f_1, f_2],id_M) \leq (1 + C(\pi(f_2^{-1})))(1 + C(\pi(f_1))) \widehat{d}_\alpha(\phi, id_M).
\end{equation}
If $\widehat{d}_\alpha(\phi, id_M) = 0$, then it follows as above that $\widetilde{d}_\alpha$ vanishes identically.\\
\end{rmk}

\subsection{Proof of the first part of Theorem \ref{thm:non-degeneracy of d_alpha non-universal}} \hfill 

The first part of Theorem \ref{thm:non-degeneracy of d_alpha non-universal} could be proven in the same way since also $\textnormal{Cont}^1_0(M,\xi)$ is perfect, but $\textnormal{Cont}^1_0(M,\xi)$ is, in fact, simple \cite{tsu08}, so the claim follows directly from the observation that $\{\phi \in \textnormal{Cont}^1_0(M,\xi)|\, d_\alpha(Id_M,\phi) = 0\}$ is a normal subgroup of $\textnormal{Cont}^1_0(M,\xi)$ by the same density argument that was used above. \eproof\\

\begin{rmk}
	A similar argument of first using $C^1$-contactomorphism and then approximating by smooth ones was used in \cite{fpr18}. As was pointed out in \cite{fpr18}, this can be avoided by using the results in \cite{ryb10} about smooth contactomorphisms instead of the case of finite regularity in \cite{tsu08}.\\
\end{rmk}

\subsection{Proof of the second parts of Theorems \ref{thm:non-deg of d_alpha} and \ref{thm:non-degeneracy of d_alpha non-universal}} \hfill

To prove the last part of Theorem \ref{thm:non-deg of d_alpha}, assume that $M$ is compact. 

Denote the Reeb flow on $(M,\alpha)$ by $\phi^\alpha_t$ and its lift to $\tildecont(M,\xi)$ by $f^\alpha_t$. If $\widetilde{d}_\alpha$ vanishes identically on $\tildecont(M,\xi)$, then there exists a Hamiltonian $H_t:M \to (-1, \infty)$ so that $[\phi^H_t] = f^\alpha_{-1}$, where $\phi^H_t$ denotes the contact isotopy associated to $H$. Then the contractible loop of contactomorphisms $\phi^\alpha_{t} \circ \phi^H_t$ is generated by the Hamiltonian $H_t + 1$ which is positive.

Conversely, let $\phi_t:M \to M, t \in [0,1],$ be a contractible positive loop of contactomorphisms. Without loss of generality, we may assume that $\phi_0 \neq \phi_{\frac{1}{2}}$. Then $\widetilde{d}_\alpha([\phi_0],[\phi_{\frac{t}{2}}]) = 0$ since there exist both positive and negative contact isotopies from $[\phi_0]$ to $[\phi_{\frac{t}{2}}]$. By what we already proved above, this implies that $\widetilde{d}_\alpha$ vanishes identically. \eproof\\

The proof of the second part of Theorem \ref{thm:non-degeneracy of d_alpha non-universal} works in the same way, except that we do not keep track of the homotopy classes of contact isotopies and positive loops.
\eproof\\

\section{Proof of Theorem \ref{thm:chekanov's dichotomy for d_alpha}}\label{sec:proof of chekanov's dichotomoy for d_alpha}

In this section, we prove Theorem \ref{thm:chekanov's dichotomy for d_alpha}. For the most part, the proof follows the arguments in \cite{rz18} and is an adaption of Appendix D of the author's master's thesis \cite{nak19} to the present setting. In fact, the degeneracy of the unparametrized metric $d^A_\alpha$ for non-coisotropic submanifolds and the vanishing of the parametrized and unparametrized metrics follow immediately from the corresponding result for the Shelukhin-Chekanov-Hofer metric that were proven in \cite{rz18} and \cite{drs24b}, and the dichotomy for unparametrized Legendrians is a special case of Lemma 7 in \cite{can23}, but since the literature does not cover all of the different cases we consider here, we include a proof for completeness. 

Let $i:A \hookrightarrow M$ be the inclusion of a closed subset. We write 
\begin{equation}
\Sigma_A \coloneqq \{\phi \in \Cont(M,\xi)|\, \phi(A) = A\}
\end{equation}
and 
\begin{equation}
\Sigma_i \coloneqq \{\phi \in \Cont(M,\xi)|\, \phi \circ i = i\}
\end{equation}
for the stabilizer subgroups of $A$ and $i$.\\

\begin{prop}\label{sigmabar}
	The closure of $\Sigma_A$ with respect to the topology on $\Cont(M,\xi)$ induced by $d_\alpha$ is given by 
	\begin{equation}
	\overline{\Sigma}_A = \{\phi \in \Cont(M,\xi)|\, d^A_\alpha(A,\phi(A)) = 0\}.
	\end{equation}
	
	The closure of $\Sigma_i$ with respect to the topology on $\Cont(M,\xi)$ induced by $d_\alpha$ is given by 
	\begin{equation}
	\overline{\Sigma}_i = \{\phi \in \Cont(M,\xi)|\, d^i_\alpha(i,\phi \circ i) = 0\}.
	\end{equation}
\end{prop}

\bproof
Let $\phi \in \overline{\Sigma}_A$ and let $\eps > 0$ be a number . Then there exists $\psi \in \Sigma_A$ with $d_\alpha(\phi \psi^{-1},id_M) = d_\alpha(\phi, \psi) < \eps$. Since $\psi(A) = A$ and thus also $\psi^{-1}(A) = A$, we see that 
\begin{equation}
d_\alpha^A(A,\phi(A)) = d_\alpha^A(A,\phi(\psi^{-1}(A))) < \eps.
\end{equation}
Because $\eps > 0$ was arbitrary, it follows that $d^A_\alpha(A,\phi(A)) = 0$.

Conversely, let $\phi \in \Cont(M,\xi)$ be a contactomorphism with $d^A_\alpha(A,\phi(A)) = 0$, i.e. given any number $\eps > 0$ there exists a contactomorphism $\psi \in \Cont(M,\xi)$ with $\psi(A) = \phi(A)$ and $d_\alpha(id_M,\psi) < \eps$. Thus, $\psi^{-1} \phi \in \Sigma_A$ and 
\begin{equation}
d_\alpha(\psi^{-1} \phi, \phi) = d_\alpha(\psi^{-1},id_M) = d_\alpha(\psi,id_M) < \eps.
\end{equation}
Since $\eps > 0$ was arbitrary, it follows that $\phi \in \overline{\Sigma}_A$.

The proof of the second part of the proposition is only notationally different.
\eproof\\

In particular, this implies that $\overline{\Sigma}_A$ and $\overline{\Sigma}_i$ are groups as the closure of any subgroup of a topological group is again a subgroup.

We define the rigid loci of $A$ and $i$ by
\begin{equation}
R_A \coloneqq \bigcap\limits_{\phi \in \overline{\Sigma}_A} \phi^{-1}(A)
\end{equation}
and 
\begin{equation}
R_i \coloneqq \bigcap\limits_{\phi \in \overline{\Sigma}_i} \phi^{-1}(A),
\end{equation}
which both are closed subsets of $A$.

For any open subset $U \subseteq M$, $\text{Cont}_U(M,\xi)$ denotes the subgroup of all contactomorphisms $\phi \in \Cont(M,\xi)$ that are generated by a Hamiltonian with support in $U$.

\begin{prop}\label{prop:cont m setminus r_a subset sigmabar}
	We have 
	\begin{equation}
	\mathrm{Cont}_{M\setminus R_A}(M,\xi) \subseteq \overline{\Sigma}_A
	\end{equation}
	and 
	\begin{equation}
	\mathrm{Cont}_{M\setminus R_i}(M,\xi) \subseteq \overline{\Sigma}_i
	\end{equation}
\end{prop}

\bproof 
Let $x \in M \setminus R_A$ be a point. By definition of $R_A$ there exists a contactomorphism $\phi \in \overline{\Sigma}_A$ with $\phi(x) \not\in A$. Since $A$ is closed, there exists a neighbourhood $U_x$ of $x$ with $\phi(U_x) \cap A = \emptyset$. Then it follows that
\begin{equation}
\phi\, \text{Cont}_{U_x}(M,\xi) \phi^{-1} = \text{Cont}_{\phi(U_x)}(M,\xi) \subseteq \overline{\Sigma}_A.
\end{equation}
Because $\overline{\Sigma}_A$ is a subgroup of $\Cont(M,\xi)$ and $\phi \in \overline{\Sigma}_A$, we see that $\text{Cont}_{U_x}(M,\xi) \subseteq \overline{\Sigma}_A$. Note that the sets $U_x,~x\in M \setminus R_A,$ cover $M \setminus R_A$. By the fragmentation lemma (see \cite{ban97}), $\text{Cont}_{M \setminus R_A}(M,\xi)$ is generated by $\bigcup\limits_{x \in M \setminus R_A} \text{Cont}_{U_x}(M,\xi)$ which proves that $\text{Cont}_{M \setminus R_A}(M,\xi) \subseteq \overline{\Sigma}_A$.

The proof of the second inclusion is only notationally different.
\eproof\\

\begin{rmk}\label{rmk:R_A=empty,R_A=A}
	Proposition \ref{prop:cont m setminus r_a subset sigmabar} implies that if $R_A = \emptyset$ (resp. $R_i = \emptyset$), then $\overline{\Sigma}_A = \Cont(M,\xi)$ (resp. $\overline{\Sigma}_i = \Cont(M,\xi)$) which means that $d^A_\alpha$ (resp. $d^i_\alpha$) vanishes identically by Proposition \ref{sigmabar}.
	
	On the other hand, it is easy to check that $R_A = A$ if and only if $d_\alpha$ is non-degenerate on $\mathscr{L}(A)$. Indeed, we see that
	\begin{equation}
	\begin{aligned}
	R_A = A &\Leftrightarrow \forall \phi \in \overline{\Sigma}_A: A \subseteq \phi^{-1}(A) \quad (\Leftrightarrow \phi(A) \subseteq A)\\
	&\Leftrightarrow  \forall \phi \in \overline{\Sigma}_A: A = \phi(A) \\
	&\Leftrightarrow d^A_\alpha \text{ is non-degenerate on } \mathscr{L}(A),
	\end{aligned}
	\end{equation}
	where we used that for any $\phi \in \overline{\Sigma}_A$ also $\phi^{-1} \in \overline{\Sigma}_A$ in order to get the equivalence between the first and the second line. The equivalence between the second and the third line follows from  Proposition \ref{sigmabar}.
	
	Similarly, 
	\begin{equation}
	\begin{aligned}
	R_i = A &\Leftrightarrow \phi \in \overline{\Sigma}_i: A \subseteq \phi^{-1}(A)\\
	&\Leftrightarrow  \forall \phi \in \overline{\Sigma}_i: A = \phi(A) \\
		&\Leftrightarrow \forall j \in \mathscr{L}(i), j(A) \neq A: d_\alpha^i(i,j) > 0\\
	&\Leftrightarrow \forall i_0,i_1 \in \mathscr{L}(i), i_0(A) \neq i_1(A): d_\alpha^i(i_0,i_1) > 0. 
	\end{aligned}
	\end{equation}\\
\end{rmk}

Proposition \ref{prop:cont m setminus r_a subset sigmabar} shows that a contactomorphism in $\Cont(M,\xi)$ that is equal to the identity in a neighbourhood of $R_A$ (resp. $R_i$) lies in $\overline{\Sigma}_A$ (resp. $\overline{\Sigma}_i$). It turns out that for this to be true we actually only have to demand that the contactomorphism is generated by a contact Hamiltonian that vanishes along $R_A$ (resp. $R_i$), as the following proposition shows:

\begin{prop}\label{prop:h vanishes on r_a}
	Let $H \in C^\infty_c(M \times [0,1])$ be a compactly supported Hamiltonian. If $H|_{R_A \times [0,1]} \equiv 0$, then its contact Hamiltonian flow $\phi^H_t$ satisfies $\phi^H_s \in \overline{\Sigma}_A$ for all $s \in [0,1]$.
	
	If $H|_{R_i \times [0,1]} \equiv 0$, then $\phi^H_s \in \overline{\Sigma}_i$ for all $s \in [0,1]$.
\end{prop}

\bproof
If we define $H^s(t,x) \coloneqq sH(st,x)$, then $\phi^{H_s}_1 = \phi^H_s$. Therefore, it will be enough to prove the proposition in the case that $s = 1$.

Let $H_t:M \to \R$ be a function with $H_t|_{R_A} \equiv 0$. For any number $n \in \N$, let $f_n:\R \to \R$ be a non-decreasing function such that $f_n(x) = x$ if $|x| \geq \frac{2}{n}$ and $f_n(x) = 0$ if $|x| \leq \frac{1}{n}$. Then $\phi^{f_n \circ H}_1$ converges to $\phi^H_1$ in the topology defined by $d_\alpha$ since $\phi_{-t}^H\phi^{f_n \circ H}_t$ is generated by the Hamiltonian $\chi_t(f_n \circ H_t - H_t) \circ \phi^H_t$, where $\chi_t$ is the positive function defined by $(\phi^H_t)^*\alpha = \chi_t \alpha$. Because $H_t$ vanishes on $R_A$ for all $t \in [0,1]$, $f_n \circ H_t$ vanishes in a neighbourhood of $R_A$ and hence,
\begin{equation}
\phi^{f_n \circ H}_1 \in \text{Cont}_{M \setminus R_A}(M,\xi) \subseteq \overline{\Sigma}_A
\end{equation}
by Proposition \ref{prop:cont m setminus r_a subset sigmabar}. As $\overline{\Sigma}_A$ is closed, it follows that also $\phi^H_1 \in \overline{\Sigma}_A$.

The proof of the second part is the same.
\eproof\\

In the following, we will only consider autonomous contact Hamiltonians $H \in C^\infty_c(M)$. For any closed subset $B \subseteq M$ we define
\begin{equation}
I_B \coloneqq \{H \in C^\infty_c(M)|\, H|_B = 0\}.
\end{equation}

Let $\{ \cdot,\cdot\}$ denote the Lie bracket defined by 
\begin{equation}\label{eq:def of lie bracket}
\begin{split}
\{F,G\} &\coloneqq \alpha([X_F,X_G]) \\
&= (dG)(X_F) -(dF)(X_G) - d \alpha(X_F,X_G) \\
&= (dG)(X_F) - (dF)(R_\alpha) G.
\end{split}
\end{equation}
for $F,G \in C_c^\infty(M)$, where $X_F$ and $X_G$ denote the contact vector fields associated to $F$ and $G$, respectively.\\

\begin{cor}\label{cor:I_R_A closed under lie bracket}
	$I_{R_A}, I_{R_i} \subseteq C_c^\infty(M)$ are closed under $\{ \cdot,\cdot\}$.
\end{cor}

\bproof Let $F,G \in I_{R_A}$. It follows immediately from the definition of $R_A$ that $R_A$ is preserved by all elements of $\overline{\Sigma}_A$. By Proposition \ref{prop:h vanishes on r_a} we know that $\phi^F_t \in \overline{\Sigma}_A$ and thus $\phi^F_t(R_A) = R_A$ for all $t \in [0,1]$. Because $G \in I_{R_A}$, this implies that $G \circ \phi^F_t|_{R_A} \equiv 0$ for all $t \in [0,1]$. Hence, for any $x \in R_A$ we have
\begin{equation}
\{F,G\}(x) = (dG)(X_F)(x) = \left.\frac{d}{dt}\right|_{t=0} (G \circ \phi^F_t)(x) = 0.
\end{equation}

The proof for $I_{R_i}$ is the same.
\eproof\\

As the following lemma shows, Corollary \ref{cor:I_R_A closed under lie bracket} restricts the possible structures of the rigid locus.

\begin{lem}\label{lem:submanifold closed under lie bracket}\hspace{-1mm}\textnormal{\cite{rz18}}\,
	If $N \subseteq M$ is a submanifold and $x \in N$ is a point, then there exist functions $F,G \in I_N$ with $\{F,G\}(x) \neq 0$ if and only if
	\begin{equation}
	(T_xN \cap \xi_x)^{\perp_{d \alpha}} \not\subseteq T_xN \cap \xi_x
	\end{equation}
	is satisfied. \eproof\\
\end{lem}

\begin{rmk}\label{rmk:n-dim coisotropic is legendrian}
	If $\dim(M) = 2n+1$ and $A$ is a submanifold of $M$ of dimension $\dim(A) = n$, then the condition in the lemma is equivalent to $A$ being non-Legendrian at $x$.\\
\end{rmk}

As an immediate consequence of Corollary \ref{cor:I_R_A closed under lie bracket} and Lemma \ref{lem:submanifold closed under lie bracket} we obtain the following.

\begin{cor}\label{cor:R_A only union of components if a coisotropic}
	If $i: A \hookrightarrow M$ is the inclusion of a properly embedded submanifold and $R_A$ (resp. $R_i$) is a union of components of $A$, then 
	\begin{equation}
	(T_xA \cap \xi_x)^{\perp_{d \alpha}} \subseteq T_xA \cap \xi_x 
	\end{equation} 
	for all $x \in R_A$ (resp. for all $x \in R_i$). \eproof \\
\end{cor}

Combined with Remark \ref{rmk:R_A=empty,R_A=A}, this implies the first assertion of Theorem \ref{thm:chekanov's dichotomy for d_alpha}:

\begin{cor}\label{cor:d only non-degenerate if a coisotropic}
	If $i: A \hookrightarrow M$ is the inclusion of a properly embedded submanifold and $d^A_\alpha$ is non-degenerate or for any $i_0,i_1 \in \mathscr{L}(i)$ with $i_0(A) \neq i_1(A)$ we have that  $d_\alpha^i(i_0,i_1) > 0$, then 
	\begin{equation}
	(T_xA \cap \xi_x)^{\perp_{d \alpha}} \subseteq T_xA \cap \xi_x 
	\end{equation} 
	for all $x \in A$. \eproof\\
\end{cor}

\begin{rmk}\label{rmk:def cosiotropic in contact manifold}
	Submanifolds satisfying the condition in the corollary are also called \emph{coisotropic} (see \cite{hua15}). \\ 
\end{rmk}

\begin{cor}\label{cor:embedding ball in r_a}
	Let $x \in R_A$ be a point and let $F_1,...,F_k \in I_{R_A}$ (resp. $F_1,...,F_k \in I_{R_i}$) be functions such that $dF_1|_{\xi_x},...,dF_k|_{\xi_x} \in \xi^*_x$ are linearly independent. Then the map
	\begin{equation}
	\begin{aligned}
	\psi: \R^k &\to M \\
	(a_1,...,a_k) &\mapsto \phi^{\sum_i a_iF_i}_1(x)
	\end{aligned}
	\end{equation}
	has image contained in $R_A$ (resp. $R_i$) and its restriction to a sufficiently small ball around the origin is an embedding.
\end{cor}

\bproof
Because $\sum_i a_i F_i \in I_{R_A}$, we see that $\psi$ has image contained in $R_A$ (cf. the proof of Corollary \ref{cor:I_R_A closed under lie bracket}).

The differential of $\psi$ at $0 \in \R^k$ maps the standard basis vectors $e_1,...,e_k$ to the vectors $X_{F_1}(x),...,X_{F_k}(x)$. The assumption on the differentials $dF_i|_{\xi_x}$ implies that these vectors are linearly independent. Thus, the restriction of $\psi$ to a sufficiently small neighbourhood of the origin is an embedding.

The proof for $R_i$ is the same.
\eproof\\

\begin{cor}\label{cor:r_a contained in l}
	Let $i:A \hookrightarrow M^{2n+1}$ be the inclusion of a closed subset of $M$ with $R_A \neq \emptyset$ (resp. $R_i \neq \emptyset$) and let $L \subseteq M$ be a properly embedded submanifold. If $\dim(L) < n$, then $R_A$ (resp. $R_i$) is not contained in $L$ and if $\dim(L) =n$ and $R_A$ is contained in $L$, then $R_A$ (resp. $R_i$) is a non-empty union of connected components of $L$.
\end{cor}

\bproof
Suppose that $R_A \subseteq L$ and denote $k \coloneqq \dim(M) - \dim(L)$. Let $x \in R_A$ be a point. Using a coordinate chart around $x$ in which $L$ looks like $\{0\} \times \R^{\dim(M)-k}$ we can find functions $F_1,...,F_k$ that vanish on $L$ such that $dF_1(x),...,dF_k(x)\in T_x^*M$ are linearly independent. This means that the linear map
\begin{equation}
(dF_1(x),...,dF_k(x)): T_xM \to \R^k 
\end{equation}
has rank $k$ which implies that the linear map 
\begin{equation}
(dF_1,...,dF_k)|_{\xi_x}:\xi_x \to \R^k
\end{equation}
has at least rank $k-1$. After reordering the $F_i$ we may assume that the linear map
\begin{equation}
(dF_1,...,dF_{k-1})|_{\xi_x}:\xi_x \to \R^k
\end{equation}
has rank $k-1$. Corollary \ref{cor:embedding ball in r_a} now implies that we can embed a $(k-1)$-dimensional ball into $R_A \subseteq L$. 

In the case that $\dim(L) < n$ this immediately gives a contradiction since 
\begin{equation}
k-1 = 2n +1 - \dim(L) -1  > n > \dim(L).
\end{equation}

If $\dim(L) = n$, then $k-1 = \dim(L)$ and we have shown that a ball-shaped neighbourhood of $x \in L$ is contained in $R_A$. This means that $R_A \subseteq L$ is open and since it is also closed and non-empty, the corollary follows.

The proof for $R_i$ is the same.
\eproof\\

Corollary \ref{cor:d only non-degenerate if a coisotropic} and Remark \ref{rmk:n-dim coisotropic is legendrian} already show that $d^A_\alpha$ is degenerate and that $d^i_\alpha$ is ``degenerate on embeddings with different images" if $A \subseteq M^{2n+1}$ is a non-Legendrian submanifold of dimension $n$. But an even stronger statement holds, which gives the second and third assertion of Theorem \ref{thm:chekanov's dichotomy for d_alpha}:  

\begin{cor}\label{cor:n-dim submafds either rigid not completely non-rigid}
	Let $i:A \subseteq M^{2n+1}$ be the inclusion of a properly embedded submanifold of dimension $\leq n$. If $A$ has no Legendrian component, then $d_\alpha^A$ and $d_\alpha^i$ vanish identically. If $A$ is a connected Legendrian submanifold, then $d^A_\alpha$ either vanishes identically or is non-degenerate, and $d_\alpha^i$ either vanishes identically or for any $i_0,i_1 \in \mathscr{L}(i)$ with $i_0(A) \neq i_1(A)$ we have that  $d_\alpha^i(i_0,i_1) > 0$.
\end{cor}

\bproof
Corollary \ref{cor:r_a contained in l} implies that $R_A$ is empty if $\dim(A) < n$ and that $R_A$ is a union of components of $A$ if $\dim(A) = n$. In the later case we know by Corollary \ref{cor:R_A only union of components if a coisotropic} (and Remark \ref{rmk:n-dim coisotropic is legendrian}) that every component of $R_A$ is Legendrian. Then Remark \ref{rmk:R_A=empty,R_A=A} finishes the proof.

The proof for $d_\alpha^i$ is analogous. 
\eproof\\

The next proposition finishes the proof of Theorem \ref{thm:chekanov's dichotomy for d_alpha}. 

\begin{prop}
	Let $i:L \hookrightarrow M$ be the inclusion of a compact and connected Legendrian submanifold. Then $d_\alpha^L$ vanishes identically if and only if there exists a positive loop of Legendrians starting at $L$. $d^i_\alpha$ vanishes identically if and only if there exists a positive loop of parametrized Legendrians starting at $i$.
\end{prop}

The proof differs only notationally from the proof of the second part of Theorem \ref{thm:non-deg of d_alpha}.\\

\begin{rmk}\label{rmk:proof of chekanovs dichotomy on universal covers}
On $\widetilde{L}(A)$ and $\widetilde{L}(i)$ one can similarly define
\begin{equation}
\begin{split}
\widetilde{\Sigma}_{A} &\coloneqq \{f \in \tildecont(M,\xi)|\, f([A]) = [A]\},\\
\widetilde{\Sigma}_i &\coloneqq \{f \in \tildecont(M,\xi)|\, f([i]) = [i]\}
\end{split}
\end{equation}
and 
\begin{equation}
\begin{split}
\widetilde{R}_A &\coloneqq \bigcap\limits_{f \in \overline{\widetilde{\Sigma}}_A} (\pi(f))^{-1}(A), \\
\widetilde{R}_i &\coloneqq \bigcap\limits_{f \in \overline{\widetilde{\Sigma}}_i} (\pi(f))^{-1}(A),
\end{split}
\end{equation}
where $[A]$ and $[i]$ denote the elements of $\widetilde{L}(A)$ and $\widetilde{L}(i)$ represented by the constant paths at $[A]$ and $[i]$, and $\pi: \tildecont(M,\xi) \to \Cont(M,\xi)$ is the projection. The proof of Proposition \ref{sigmabar} still applies to show that
\begin{equation}\label{eq:closure of sigmatilde A}
	\overline{\widetilde{\Sigma}}_A = \{f \in \tildecont(M,\xi)|\, \widetilde{d}^A_\alpha([A],f([A])) = 0\}
\end{equation}
and
\begin{equation}\label{eq:closure of sigmatilde i}
	\overline{\widetilde{\Sigma}}_i = \{f \in \tildecont(M,\xi)|\, \widetilde{d}^i_\alpha([i],f([i])) = 0\}.
\end{equation}
	
Assume that $i:A \to M$ is the proper embedding of a submanifold. Then the above arguments apply to show that $\widetilde{R}_A$ and $\widetilde{R}_i$ are empty if $\dim(A) < n$, that $\widetilde{R}_A = A$ or $\widetilde{R}_i = A$ only if $A$ is coisotropic, and that $\widetilde{R}_A$ and $\widetilde{R}_i$ are unions of Legendrian components of $A$ if $\dim(A) = n$.

As before, we have that $\widetilde{d}^A_\alpha$ vanishes identically if $\widetilde{R}_A = \emptyset$ and $\widetilde{d}^i_\alpha$ vanishes identically if $\widetilde{R}_i = \emptyset$.

On the other hand, observe that
\begin{equation}
	\begin{aligned}
	\widetilde{R}_A = A &\Leftrightarrow \forall f \in \overline{\widetilde{\Sigma}}_A: A \subseteq (\pi(f))^{-1}(A) \\
	&\Leftrightarrow \forall f \in \overline{\widetilde{\Sigma}}_A: A = (\pi(f))(A) \\
	&\Leftrightarrow \widetilde{d}^A_\alpha \text{ is } \pi\text{-non-degenerate on } \widetilde{\mathscr{L}}(A),
	\end{aligned}
	\end{equation}
	where the equivalence between the second and the third line follows from (\ref{eq:closure of sigmatilde A}).
	
	Similarly, 
	\begin{equation}
	\begin{aligned}
	\widetilde{R}_i = A &\Leftrightarrow \forall f \in \overline{\widetilde{\Sigma}}_i: A \subseteq (\pi(f))^{-1}(A)\\
	&\Leftrightarrow  \forall f \in \overline{\widetilde{\Sigma}}_i: A = (\pi(f))(A) \\
		&\Leftrightarrow \forall \iota \in \mathscr{\widetilde{L}}(i), (\pi(\iota))(A) \neq A: \widetilde{d}_\alpha^i([i],\iota) > 0\\
	&\Leftrightarrow \forall \iota_0,\iota_1 \in \mathscr{\widetilde{L}}(i), (\pi(\iota_0))(A) \neq (\pi(\iota_1))(A): \widetilde{d}_\alpha^i(\iota_0,\iota_1) > 0,
	\end{aligned}
	\end{equation}
	where here, $\pi$ denotes the projection $\pi: \mathscr{\widetilde{L}}(i) \to \mathscr{L}(i)$, and the equivalence between the second and the third line follows from (\ref{eq:closure of sigmatilde i}).\\
\end{rmk}

\bibliography{references}
\bibliographystyle{amsalpha}

\end{document}